\documentclass[11pt,reqno]{amsart}

\usepackage{amsthm}
\usepackage{amssymb}
\usepackage{latexsym}
\usepackage{multicol}
\usepackage{verbatim,enumerate}
\usepackage{hyperref}
\usepackage{nicematrix}
\usepackage{amsmath, amscd}
\usepackage{mathrsfs}
\usepackage[all,cmtip]{xy}
\usepackage{soul}

\advance\textwidth by 1.2in \advance\oddsidemargin by -.6in \advance\evensidemargin by -.6in
\usepackage{tikz}
\usetikzlibrary{decorations.markings}
\tikzstyle{vertex}=[ circle, fill, draw, inner sep=0pt, minimum size=4pt,]
\tikzstyle{edge}= [thick]

\newtheorem*{cor}{Corollary}
\newtheorem*{lem}{Lemma}
\newtheorem*{prop}{Proposition}

\theoremstyle{definition} \newtheorem*{defn}{Definition}
\theoremstyle{definition}
\newtheorem{thm}{Theorem}
\newtheorem*{thm*}{Theorem}

\newtheorem*{rem}{Remark}

\newenvironment{pf}{\proof}{\endproof}
\newcounter{cnt}
\newenvironment{enumerit}{\begin{list}{{\hfill\rm(\roman{cnt})\hfill}}{%
\settowidth{\labelwidth}{{\rm(iv)}}\leftmargin=\labelwidth%
\advance\leftmargin by \labelsep\rightmargin=0pt\usecounter{cnt}}}{\end{list}} \makeatletter
\def\mydggeometry{\makeatletter\dg@YGRID=1\dg@XGRID=20\unitlength=0.003pt\makeatother}
\makeatother \theoremstyle{remark}

\numberwithin{equation}{section}

\newcommand{\id}{\operatorname{id}}

\newcommand{\sgn}{\operatorname{sgn}}
\newcommand{\wt}{\operatorname{wt}}

\newcommand{\nc}{\newcommand}
\newcommand{\rnc}{\renewcommand}

\nc{\cal}{\mathcal} \nc{\goth}{\mathfrak} \rnc{\bold}{\mathbf}

\newcommand{\supp}{\operatorname{supp}}

\renewcommand{\Bbb}{\mathbb}
\nc\bomega{{\mbox{\boldmath $\omega$}}} \nc\bpsi{{\mbox{\boldmath $\Psi$}}}
 \nc\sing{{\rm sing}}
 \nc\balpha{{\mbox{\boldmath $\alpha$}}}
 \nc\bbeta{{\mbox{\boldmath $\beta$}}}
 \nc\bpi{{\mbox{\boldmath $\pi$}}}
  \nc\bpis{{\mbox{\boldmath \scriptsize$\pi$}}}
 \nc\bullets{{\mbox{\scriptsize $\bullet$}}}
 \nc\bvarpis{{\mbox{\boldmath \scriptsize$\varpi$}}}
  \nc\bvarpi{{\mbox{\boldmath $\varpi$}}}

\nc\bepsilon{{\mbox{\boldmath $\epsilon$}}}

  \nc\bomegas{{\mbox{\boldmath\scriptsize $\omega$}}}
  \nc\bepsilons{{\mbox{\boldmath \scriptsize$\epsilon$}}}
\nc{\spi}{{\rm sp}}
\nc\hlien{\hat{\lie n}^+}
  \nc\btaus{{\mbox{\boldmath \scriptsize$\tau$}}}\nc\bxi{{\mbox{\boldmath $\xi$}}}
\nc\bmu{{\mbox{\boldmath $\mu$}}} \nc\bcN{{\mbox{\boldmath $\cal{N}$}}} \nc\bcm{{\mbox{\boldmath $\cal{M}$}}} \nc\blambda{{\mbox{\boldmath
$\lambda$}}}
\nc\btau{{\mbox{\boldmath
$\tau$}}}

\newcommand{\lie}[1]{\mathfrak{#1}}

\makeatletter
\def\section{\def\@secnumfont{\mdseries}\@startsection{section}{1}%
  \z@{.7\linespacing\@plus\linespacing}{.5\linespacing}%
  {\normalfont\scshape\centering}}
\def\subsection{\def\@secnumfont{\bfseries}\@startsection{subsection}{2}%
  {\parindent}{.5\linespacing\@plus.7\linespacing}{-.5em}%
  {\normalfont\bfseries}}
\makeatother

 \nc{\Hom}{\operatorname{Hom}}
  \nc{\mode}{\operatorname{mod}}
\nc{\End}{\operatorname{End}} \nc{\wh}[1]{\widehat{#1}} \nc{\Ext}{\operatorname{Ext}}
 \nc{\ch}{\operatorname{ch}} \nc{\ev}{\operatorname{ev}}
\nc{\Ob}{\operatorname{Ob}} \nc{\soc}{\operatorname{soc}} \nc{\rad}{\operatorname{rad}} \nc{\head}{\operatorname{head}}

 \nc{\Cal}{\cal} \nc{\Xp}[1]{X^+(#1)} \nc{\Xm}[1]{X^-(#1)}
\nc{\on}{\operatorname} \nc{\Z}{{\bold Z}} \nc{\J}{{\cal J}} \nc{\C}{{\bold C}} \nc{\Q}{{\bold Q}}

\nc{\N}{{\Bbb N}} \nc\boa{\bold a} \nc\bob{\bold b} \nc\boc{\bold c} \nc\bod{\bold d} \nc\boe{\bold e} \nc\bof{\bold f} \nc\bog{\bold g}
\nc\boh{\bold h} \nc\boi{\bold i} \nc\boj{\bold j} \nc\bok{\bold k} \nc\bol{\bold l} \nc\bom{\bold m} \nc\bon{\bold n} \nc\boo{\bold o}
\nc\bop{\bold p} \nc\boq{\bold q} \nc\bor{\bold r} \nc\bos{\bold s} \nc\boT{\bold t} \nc\boF{\bold F} \nc\bou{\bold u} \nc\bov{\bold v}
\nc\bow{\bold w} \nc\boz{\bold z} \nc\boy{\bold y} \nc\ba{\bold A} \nc\bb{\bold B} \nc\bc{\mathbb C} \nc\bd{\bold D} \nc\be{\bold E} \nc\bg{\bold
G} \nc\bh{\bold H} \nc\bi{\bold I} \nc\bj{\bold J} \nc\bk{\bold K} \nc\bl{\bold L} \nc\bm{\bold M}  \nc\bo{\bold O} \nc\bp{\bold
P} \nc\bq{\bold Q} \nc\br{\bold R} \nc\bs{\bold S} \nc\bt{\bold T} \nc\bu{\bold U} \nc\bv{\bold V} \nc\bw{\bold W} \nc\bx{\bold
x} \nc\KR{\bold{KR}} \nc\rk{\bold{rk}} \nc\het{\text{ht }}
\nc\bz{\mathbb Z}
\nc\bn{\mathbb N}
\nc\us{\underline \bos}
\nc\uS{ \bs_{{\rm alt}}}
\nc\pr{\rm pr}

\nc\toa{\tilde a} \nc\tob{\tilde b} \nc\toc{\tilde c} \nc\tod{\tilde d} \nc\toe{\tilde e} \nc\tof{\tilde f} \nc\tog{\tilde g} \nc\toh{\tilde h}
\nc\toi{\tilde i} \nc\toj{\tilde j} \nc\tok{\tilde k} \nc\tol{\tilde l} \nc\tom{\tilde m} \nc\ton{\tilde n} \nc\too{\tilde o} \nc\toq{\tilde q}
\nc\tor{\tilde r} \nc\tos{\tilde s} \nc\toT{\tilde t} \nc\tou{\tilde u} \nc\tov{\tilde v} \nc\tow{\tilde w} \nc\toz{\tilde z} \nc\woi{w_{\omega_i}}
\nc\chara{\operatorname{Char}}

\begin{document}

\title[]{Alternating snake modules and a determinantal formula}
\author{Matheus Brito}
\address{Departamento de Matematica, UFPR, Curitiba - PR - Brazil, 81530-015}
\email{mbrito@ufpr.br}
\thanks{M.B. was partially supported CNPq grant 405793/2023-5}
\author{Vyjayanthi Chari}
\address{Department of Mathematics, University of California, Riverside, 900 University Ave., Riverside, CA 92521, USA}
\email{chari@math.ucr.edu}
\thanks{V.C. was partially supported by a travel grant from the Simons Foundation and by the Infosys Visiting Chair position at the Indian Institute of Science, Bangalore.}

\begin{abstract} We introduce a family of modules for the quantum affine algebra which include as very special cases both the snake modules  and  modules arising from a  monoidal categorification of cluster algebras. We give necessary and sufficient conditions for these modules to be prime and prove a unique factorization result. We also give an explicit formula expressing the module as an alternating sum of Weyl modules. Finally, we give an application of our results to a classical question in  the category $\cal O(\lie{gl}_r)$. Specifically we apply our results to show that there are a large family of non--regular, non--dominant weights $\mu$ for which the non--zero Kazhdan--Lusztig coefficients $c_{\mu, \nu}$ are $\pm 1$.
\end{abstract}

\maketitle

\section*{Introduction}
The study of finite--dimensional representations of a quantum affine algebra has been a central topic in representation theory for over three decades. The subject has  deep connections to various fields, including integrable systems, algebraic geometry, and mathematical physics. More recently the connection with cluster algebras through the work of \cite{HL10, HL13a} has brought many new ideas to the subject. The work of \cite{KKKO18, KKOP, KKOP24, KKOP2} has led to remarkable developments in the area and  new tools are now available for the study of these representations.\\\\
In their papers, Hernandez and Leclerc identified a certain 
 tensor subcategory denoted $\mathscr F_n$  of the category of  finite--dimensional representations of the quantum affine algebra. They showed that there was an isomorphism between the  Grothendieck ring of this category and an infinite rank cluster algebra. They conjectured, now a theorem \cite{KKKO18, KKOP, KKOP24, KKOP22, Qin17} that a cluster monomial corresponds to an irreducible representation whose tensor square is irreducible; such representations are called real. Moreover a cluster variable corresponded to an irreducible representation which is not isomorphic to a tensor product of nontrivial irreducible representations; such representations are called prime. They also conjectured the converse; namely all real representations in the category are cluster monomials and real prime representations are cluster variables. But this is only known to be true for  very specific families of representations and is open in general. One of the reasons for this,  is that it  is highly nontrivial to prove that a module is  prime or real. For some combinatorial approaches to the problem of classifying prime representations  see \cite{EL24, MS24}.
\\\\
From now on we restrict our attention to a quantum affine algebra of  type $A_n$. In this case the  irreducible modules in the Hernandez--Leclerc subcategory are indexed by a free abelian monoid $\cal I_n^+$ generated by elements $\bomega_{i,j}$ where $i,j\in\mathbb Z$ and $0\le j-i\le n+1$. (This is a reformulation of the usual index set: the  Drinfeld polynomials.) Associated with every element of this monoid one also has a standard or Weyl module.
An important family of real modules which are known (\cite{DLL19}) to be cluster monomials    are the snake modules introduced by Mukhin and Young.  These are indexed by elements of the form $\bomega_{i_1,j_1}\cdots\bomega_{i_r,j_r}$ with $i_1<\cdots<i_r$ and $j_1<\cdots<j_r$. These modules have many nice properties and their characters are explicitly known. 
\\\\
The index set for snake modules also defines a family of modules  in the category $\mathcal M_N$ of finite length, complex  smooth representations  of $GL_N(F)$ where $F$ is a non--Archimedean field; in that context they are called the ladder modules and have been studied in \cite{Gur21, LM14, LM18}. The irreducible representations in this world are also indexed by elements of $\cal I_n^+$; here the index set consists of the Zelevinsky multisegments. There is an associated notion of square irreducibility which is the analog of real modules in the quantum setting.
\\\\
Loosely speaking,  one can use an  affine Schur Weyl duality to go between the  category $\mathscr F_n$ and the Bernstein block in $\mathcal M_N$;  the snake modules correspond to the ladder modules. In \cite{BaC15} the authors explained the connection between $\mathcal M_N$ and the BGG--category  $\cal O$ for $\lie {gl}_r$. In particular the BGG--resolution of a finite--dimensional irreducible module of $\lie{gl}_r$ gives a resolution of the irreducible  ladder modules in terms of standard modules. Using \cite{CPHecke}  one can show that this leads to a resolution of the snake module by Weyl modules. \\\\
In \cite{LM18}, Lapid and M\'inguez continued their study of  smooth complex representations of $GL_N(F)$. They give several equivalent definitions for an irreducible representation associated to a regular element  to be square irreducible.  A regular element is an element of the form $\bomega_{i_1,j_1}\cdots\bomega_{i_r, j_r}$ where  $i_s\ne i_p$ and $j_s\ne j_p$ for all $1\le p\ne s\le r$. They show that the property of square irreducibility also holds for certain non--regular representations.\\\\
In the quantum affine setting  there are interesting representations coming from the connection with the cluster algebras \cite{BC19a, HL10, HL13}  which  are not regular. In the current paper we introduce a family of modules which we call alternating snake modules. The snake modules and the modules coming from the category $\cal C_1$ of \cite{HL10} are both  very special examples of alternating snake modules. A  straightforward application of the results of \cite{KKOP} show that the modules are real. More interestingly, we give necessary and sufficient conditions for an alternating snake module to be prime. We  prove a unique factorization result; namely that an alternating snake module is isomorphic, uniquely (up to a permutation)  to a tensor product of prime alternating snake modules. Further results include a presentation of these modules, analogous to the one given in \cite{Tad95} and later generalized in \cite{LM14} for ladder modules.\\\\
We also  prove a determinantal formula for these modules (under a mild condition). Namely we define  a matrix with entries in the commutative Grothendieck ring $\cal K_0(\mathscr F_n)$ whose determinant is an alternating sum of classes of Weyl (standard) modules and  equal to the class of the irreducible module.  Under suitable conditions on the alternating snake (but still weaker than the condition that the corresponding Zelevinsky multisegments is regular) we show that the standard modules which occur with non--zero coefficients in the determinant are $\pm 1$. \\\\
Finally we give an application to the category $\cal O(\lie {gl}_r)$; namely we are able to use our result to compute in $\cal K(\cal O(\lie {gl}_r))$ the expression for certain infinite--dimensional irreducible modules in terms of the Verma modules.\\\\
For the readers convenience, we establish in the first section, the minimal possible notation to define the notion of alternating snakes, give examples and state all the main results including the connection with $\cal O(\lie{gl}_r)$. The proofs are given in the subsequent sections.  \\\\

{\em Acknowledgment: The authors thank David Hernandez and Bernard Leclerc for many helpful discussions and insightful questions. They  thank Ryo Fujita for drawing their attention to the connection with the Arakawa--Suzuki functor and for pertinent references. 
A substantial portion of this work was carried out during two visits to  Oberwolfach as part of the OWRF program in 2023; the authors are deeply appreciative of the excellent environment at the Mathematisches Forschungsinstitut. M.B. is grateful to the Department of Mathematics, UCR, for their
hospitality during a visit when part of this research was carried out.}

\section{Alternating Snake Modules: Main Results } We begin by recalling  some essential definitions and results on the representation theory of the quantum loop algebra $\widehat\bu_n$ associated to $\lie{sl}_{n+1}$. We then introduce a new   family of irreducible modules for $\widehat\bu_n$ which we call alternating snake modules. After that we   state the main results of the paper and  end the section with an application of our results to the BGG-category $\cal O$ for the Lie algebra $\lie{gl_r}$.\\\\
 Assume throughout that $q$ is a non--zero complex number and not a root of unity. As usual $\mathbb C$ (resp. $\mathbb C^\times$, $\mathbb Z$, $\mathbb Z_+$, $\mathbb N$) will denote the set of complex numbers (resp. non-zero complex numbers, integers, non-negative integers, positive integers).  Given $\ell\in \mathbb N$ we denote by $\Sigma_\ell$ the symmetric group on $\ell$ letters. 

\subsection{The algebra $\widehat\bu_n$ and the category $\mathscr F_n$}\label{basicdef}  For $n\in\mathbb N$, let  $\widehat\bu_n$ be  the quantum loop algebra associated to $\lie{sl}_{n+1}(\mathbb C)$;  we refer the reader to \cite{CP95} for precise definitions. For our purposes, it is enough to recall that $\widehat\bu_n$ is a Hopf algebra with 
an infinite set of generators: $x_{i,s}^\pm$, $k_i^{\pm 1}$, $\phi^\pm_{i,s}$,  $1\le i\le n$ and $s\in\mathbb Z$. The subalgebra $\widehat\bu_n^0$   generated by the elements $\phi^\pm_{i,s}$, $1\le i\le n$,  $s\in\mathbb Z$ is commutative.
\\\\
It is well known (see \cite{CP91,CP95}) that the isomorphism  classes of irreducible finite--dimensional representations  of $\widehat\bu_n$  are parameterized  by elements of a free abelian monoid with identity $\bold 1$   and generators  $\bvarpi_{m,a}$ with   $1\le m\le n$ and $a\in\mathbb C^{\times}$. The trivial representation of $\widehat\bu_n$ corresponds to the identity element of the monoid. It was shown in \cite{CP01} that corresponding to an element of this monoid there also exists a finite--dimensional indecomposable module called a Weyl module which has the corresponding irreducible module as its unique irreducible quotient. 
\\\\
Let $\mathscr F_{n}$ be the full subcategory of the category of finite--dimensional representations of $\widehat\bu_n$ consisting of objects whose  Jordan--Holder components are indexed by the submonoid (with identity) generated by elements $\bvarpi_{m,q^a}$ with  $a-m\in2\mathbb Z$. It was proved in \cite{HL10} that $\mathscr F_{n}$ is a rigid tensor category and  we  let $\cal K_0(\mathscr F_{n})$ be the corresponding  Grothendieck ring. The results of \cite{FR99} show that this     ring is commutative with basis given by the classes of the simple objects.  For any object $V$ of $\mathscr F_n$ we denote by $[V]$ the corresponding element of $ \cal K_0(\mathscr F_n)$. 
\subsection{The group $\cal I_n$} It will be convenient  to use a different index set for the simple objects of $\mathscr F_n$. 
Let  $\mathbb I_n$ be the set of intervals $[i,j]$ with $i,j\in\mathbb Z$ and $0\le j-i\le n+1$ and for $r\ge 1$ let $\mathbb I_n^r$ be the set of ordered $r$--tuples of elements of $\mathbb I_n$. Given elements $\bos_1\in\mathbb I_n^{r_1}$ and $\bos_2\in\mathbb I_n^{r_2}$ we let $\bos_1\vee\bos_2$ be the element of $\mathbb I_n^{r_1+r_1}$ obtained by concatenation.\\\\
Define $\cal I_n^+$ (resp. $\cal I_n$) to  be the  free abelian monoid (resp. group) with identity $\bold 1$ and generators
 $\bomega_{i,j}$ with $[i,j]\in\mathbb I_n$. We understand that $\bomega_{i,i}=\bomega_{i,i+n+1}=\bold 1$ for all $i\in \mathbb Z$.
 We have a map $\mathbb I_n^r\to\cal I_n^+$ given by  $ \bos=([i_1,j_1],\cdots, [i_r,j_r])\mapsto \bomega_{\bos}=\bomega_{i_1,j_1}\cdots\bomega_{i_r, j_r}$.
Identifying a  pair $(m,q^a)$ with $1\le m\le n$ and   $a-m\in 2\mathbb Z$ with the interval $[\frac12(a-m), \frac12(a+m)]$ and $\bvarpi_{m,q^a}$ with $\bomega_{\frac 12(a-m), \frac12(a+m)}$  we see that the irreducible objects in $\mathscr F_{n}$ are also indexed by elements of $\cal I_n^+$.\\\\
Given $\bomega\in\cal I_n^+$ we let $W(\bomega)$ and $V(\bomega)$ be the  Weyl module (see Section \ref{weylpod} for the definition)  (up to isomorphism)  and irreducible module in $\mathscr F_n$ respectively.

\subsection{$\ell$--weights} \label{ellwtsprop} 
   It was proved in  \cite{FR99} that  an object $V$  of $\mathscr F_{n}$ is  the  direct sum of generalized eigenspaces for the $\widehat{\bu}_n^0$--action. The eigenvalues are indexed by elements of $\cal I_n$ and we have, $$V=\bigoplus_{\bomegas\in\cal I_n} V_\bomegas,\ \ \wt_\ell V=\{\bomega\in\cal I_n: V_\bomegas\ne 0\},\ \ \wt_\ell ^\pm V=\wt_\ell V\cap (\cal I^+_n)^{\pm 1}.$$
     Moreover, if  $V'$ is another  object of $\mathscr F_{n}$ then  \begin{gather}
   \label{qinv}
   [V]=[V']\implies 
   \wt_\ell V=\wt_\ell V',\ \ \dim V_{\bomegas}= \dim V'_{\bomegas},\ \ \bomega\in\cal I_n,\\
\label{lwtprod}\wt_\ell (V\otimes V')=\wt_\ell V\wt_\ell V',\ \ \dim (V\otimes V')_\bomegas=\dim  (V'\otimes V)_\bomegas,\ \ \bomega\in\cal I_n,
   \end{gather}
 
\subsection{Alternating snakes} \label{altsnakesdef}
  Set \begin{gather*}\label{snake}
  \bs=\{([i_1,j_1],\cdots,[i_r,j_r])\in\mathbb I_n^r: \ \ r\ge 1,\ \ i_1<i_2<\cdots<i_r,\ \ j_1<j_2<\cdots <j_r\},
\\
\bs^{\circ}=\{([i_1,j_1],\cdots,[i_r,j_r])\in\mathbb I_n^r: \ \ r\ge 1,\ \ ([i_r,j_r],\cdots,[i_1,j_1])\in \bs\}.
\end{gather*} 
The elements of $\bs$ were called snakes in \cite{MY12} and ladders in \cite{LM14}. 
For $\bos=([i_1,j_1],\cdots,[i_r, j_r])\in\mathbb I_n$ and $0\le p<\ell\le r$, let 
\begin{equation}\label{spl}\bos(p,\ell)=([i_{p+1},j_{p+1}],\cdots, [i_{\ell}, j_{\ell}])\in\mathbb I_n^{\ell-p}.\end{equation}
We say that the elements $[i_1, j_1]$ and $[i_2,j_2]$ of $\mathbb I_n$  overlap if for some $\epsilon\in\{0,1\}$ we have
\begin{equation}\label{overlapdef} i_{1+\epsilon}<i_{2-\epsilon}\le j_{1+\epsilon}<j_{2-\epsilon}.\end{equation} Otherwise, we say that they do not overlap.
\begin{defn}\label{altsnakea}  We say that $\bos=([i_1,j_1], \cdots, [i_r,j_r])\in\mathbb I_n^r$ is an alternating snake if the following hold:
\begin{enumerate}
\item[(i)] for $1\le s\ne p\le r$ we have either $i_s\ne i_p$ or $j_s\ne j_p$,
    \item[(ii)]  the element  $\bos(s-1,s+1)$ is in $\bs^\circ\sqcup\bs$ for all $1\le s\le r-1$,
    \item[(iii)] if $1\le s<p\le r$ is such that $\bos(s-1,p)\notin\bs^\circ\sqcup\bs$ then $[i_s,j_s]$ and $ [i_p,j_p]$ do not overlap.
\end{enumerate}\hfill\qedsymbol
\end{defn}
Let $\uS$ denote the set of alternating snakes. Clearly $\bos\in\uS$ if and only if $\bos(p,\ell)\in\uS$ for all $1\le p< \ell\le r$. The modules $V(\bomega_\bos)$ with $\bos\in\uS$ are called alternating snake modules. Given $\bos\in\uS$ we define the integer $r_1:=r_1(\bos)$ to be maximal so that $\bos(0,r_1)\in\bs^\circ\sqcup\bs$.

\subsubsection {Examples}
\begin{enumerit}
    \item[(i)] 
The element 
 $
 \bos=    ([0,4], [-1,1], [1,2],[2,3])\in\mathbb I_n^4$ is an alternating snake. Note that $$\bos(0,2)=([0,4], [-1,1])\in\bs^\circ,\ \ \bos(1,4)=([-1,1], [1,2], [2,3])\in\bs, \ \ \bos(0,m)\notin\bs^\circ\sqcup\bs, \ \ m=3,4$$ and the interval $[0,4]$ does not overlap either $[1,2]$ or $[2,3]$. \\
\item[(ii)] For $n\gg 0$ and for $p\in\bz_+$  let $$\bos=( [-p, p+1], [-p+1, p+3], [-p-1, p+2], [-p, p+4], [-p-2, p+3],\cdots )\in\mathbb I_n^r.$$ Then $\bos$ is an alternating such that 
$$\bos(2k,2k+2)\in \bs, \ \ \bos(2k+1,2k+3)\in \bs^\circ, \ \ k\geq 0, \ \ \bos(m-1,m+2)\notin \bs^\circ\sqcup\bs, \ \ 1\leq m\leq r-2.$$
\item[(iii)] Suppose that $(\mu_1,\cdots, \mu_r)\in\mathbb Z^r$ and $(\lambda_1,\cdots,\lambda_r)\in\mathbb Z^r$ satisfy the following: \begin{gather*} \mu_1\le \mu_2<\mu_3\le \mu_4<\cdots,\ \ \ \lambda_1>\lambda_2\ge\lambda_3>\lambda_4\ge\cdots,\\
n+1>\lambda_1-\mu_1\ge \lambda_r-\mu_r>0.\end{gather*}
Then    $$\bos=([\mu_1,\lambda_2],\  [\mu_3,\lambda_1], [\mu_2,\lambda_4],\cdots [\mu_{2s+1},\lambda_{2s-1}],[\mu_{2s},\lambda_{2s+2}],\cdots )$$ is an alternating snake such that $\bos(m-1,m+2)\notin \bs^\circ\sqcup\bs$, for $1\leq m\leq r-2$. \\
\end{enumerit}

Further examples of alternating snakes  can be found in Section \ref{exmore}.

\subsubsection{}    Alternating snake  modules are known to be real by the work of \cite{BMS24}. Since that work is rather abstract and the proof in our case is very brief we include it  in Section \ref{sreal}.

\subsection {Prime factorizations} \label{primedecompssec} An irreducible  module  in $\mathscr F_n$ is said to be prime if  it is not isomorphic to  a tensor product of non--trivial   representations.  Clearly any irreducible object of $\mathscr F_n$ is isomorphic to  a tensor product of prime representations. It is not known in general if such a factorization is unique. 
\\\\
Our next results  show that an alternating snake module  is isomorphic (uniquely upto a permutation) to  a tensor product of prime alternating snake modules.  It also  gives a necessary and sufficient condition for $V(\bomega_\bos)$ to be prime. 
\begin{thm}\label{primefactor}
    Suppose that $\bos = ([i_1,j_1],\cdots, [i_r,j_r]) \in\uS$, $r\ge 1$. 
    \begin{enumerit}
        \item[(i)] 
   The module
    $V(\bomega_{\bos})$ is prime if the following conditions hold:  \begin{gather}\label{conn} 0\le\min\{j_s-i_{s+1}, j_{s+1}-i_s\}\le \max\{j_s-i_{s+1}, j_{s+1}-i_s\}\le n+1,\ \ 1\le s\le r-1,\\
\label{ijneq} i_{p-1}\ne i_{p+1}\ \ {\rm and}\ \ j_{p-1}\ne j_{p+1}\ \  {\rm for }\ \ 2\le p\le r-1.
\end{gather}
\item[(ii)] Suppose that  $1\le p< r$ is such that $[i_p,j_p]$ and $[i_{p+1}, j_{p+1}]$ do not satisfy \eqref{conn}. Then
\begin{gather*}V(\bomega_{\bos})\cong V(\bomega_{\bos(0,p)})\otimes  V(\bomega_{\bos(p,r)}).
\end{gather*}
\item[(iii)] Suppose that $\bos$ satisfies \eqref{conn} and  that \eqref{ijneq} does not hold for some  $2\le p'\le r-1$. There exists  $b\in\{i,j\}$  such that  $b_{p'-1}\ne b_{p'+1}$ and if we  choose $\epsilon\in\{0,1\}$ so that 
\begin{gather}\label{pf2b}
 b_{p'-1+2\epsilon}<b_{p'+1-2\epsilon}\ \ {\rm if}\ \  \bos(p'-2, p')\in\bs^\circ ,\\
 \label{pf2c}  b_{p'+1-2\epsilon}<b_{p'-1+2\epsilon}
 \ \ {\rm if}\ \ \bos(p'-2, p')\in\bs,
\end{gather}
 then $$V(\bomega_\bos)\cong V(\bomega_{\bos(0,p'-\epsilon)})\otimes V(\bomega_{\bos(p'-\epsilon, r)}).$$ 
\end{enumerit}

In particular, $V(\bomega_\bos)$ is prime if and only if \eqref{conn} and \eqref{ijneq} hold.
\end{thm}

\subsection{Prime factors} \label{pfsection} In view of the preceding theorem, it is natural to define  a prime  alternating snake to be an element of $\uS$ which satisfies \eqref{conn} and \eqref{ijneq}. Let $\uS^{\pr}$ be the set of prime alternating snakes.  It is also convenient to say that $\bos$ is connected if  it satisfies \eqref{conn}
\\\\
In the case when $\bos\in\uS\setminus\uS^{\pr}$ the preceding theorem tells us that it is natural to define the notion of a prime factor of an alternating snake. This is made precise as follows. 
 \begin{defn} \label{pfdef} 
 We say that $\bos(0,p)$ for $1\le p\le r$ is a prime factor of $\bos$ if  $\bos(0,p)\in\uS^{\pr}$ and either 
\begin{itemize}
\item
$([i_p, j_p],[i_{p+1}, j_{p+1}])$ is not connected,
 \item or $p=p'-\epsilon$ where  $p'$ and $\epsilon\in\{0,1\}$ satisfy the conditions in Theorem \ref{primefactor}(iii).
 \end{itemize}
Writing $\bos=\bos(0,p)\vee\bos(p,r)$ the remaining prime factors of $\bos$ are defined to be the  set of prime factors of $\bos(p,r)$. Clearly the prime factors come with a canonical order and 
we  call this the prime decomposition of $\bos$.  \\\\
For $1\le \ell\le \ell'\le r$ we say that $\bos(\ell-1, \ell')$ is contained in a prime factor  of $\bos$ if there exists $1\le p\le \ell\le \ell'\le p'\le r$ such that $\bos(p-1,p')$ is a prime factor of $\bos$.  Otherwise we say that $\bos(\ell-1,\ell')$ is not contained in a prime factor  of $\bos$.\end{defn}
We have the following corollary of Theorem \ref{primefactor}.
\begin{cor}\label{tensorfactor} Suppose that $\bos=\bos^1\vee\cdots\vee\bos^\ell$ is  the prime decomposition of $\bos$.  Then
\begin{gather}\label{primedecomps}
V(\bomega_\bos)\cong V(\bomega_{\bos^1})\otimes\cdots\otimes V(\bomega_{\bos^\ell}).
\end{gather} 
 Moreover  if $V(\bomega_\bos)\cong V(\bomega_1)\otimes \cdots\otimes V(\bomega_p)$ for  prime modules $V(\bomega_\ell)$, $1\le \ell\le p$ then $p=\ell$ and $\{\bomega_1,\cdots,\bomega_\ell\}=\{\bomega_{\bos^1},\cdots,\bomega_{\bos^\ell}\}$. 
\end{cor}

\subsection{A presentation of $V(\bomega_\bos)$}\label{presthm}  
 Given  $1\le p\le r-1$ such that $\bos(p-1,p+1)$ is contained in a prime factor of $\bos$ we set, $$\tau_p\bos=\bos(0,p-1)\vee ([i_{p+1}, j_{p}],[i_{p}, j_{p+1}])\vee\bos(p+1,r).$$ 

\begin{thm}\label{pres}
    Let $\bos \in\uS$.
    \begin{enumerit}
    \item[(i)]  Suppose that $1\le p\le r-1$ is such that $\bos(p-1,p+1)$ is contained in a prime factor of $\bos$. 
    Then, \begin{gather*}
\dim\Hom_{\widehat\bu_n}(W(\bomega_{\tau_p\bos}), W(\bomega_{\bos}))=1,\end{gather*}  and any non--zero element of the space is injective.    \item[(ii)]  
For $1\le p\le r-1$ let $M_p(\bos)$ be the image of a non--zero element of $ \Hom_{\widehat\bu_n}(W(\bomega_{\tau_p\bos}), W(\bomega_{\bos}))$ if $\bos(p-1,p+1)$ is contained in a prime factor of $\bos$ 
 and otherwise let $M_p(\bos)=0$. Then,
    $$V(\bomega_\bos)\cong \frac {W(\bomega_\bos)}{\sum_{p=1}^{r-1}M_p(\bos)}.$$\end{enumerit}
\end{thm}
\begin{rem} This generalizes the result of Tadi\'c (\cite{Tad95}) and Lapid--M\'inguez (\cite{LM14}) on  ladder modules.\end{rem}

\subsection{A determinantal formula} Our final result expresses $[V(\bomega_\bos)]$ (with suitable restrictions on $\bos$) as the determinant of a matrix whose entries are either zero or the elements $[V(\bomega_{i,j})]$ for some $[i,j]\in\mathbb I_n$.

\subsubsection{{\bf The matrix $A(\bos)$}} 
 We define an $r\times r$ matrix $A(\bos)$ with coefficients in $\cal K_0(\mathscr F_n)$ by induction on $r$. If $\bos=([i_1,j_1])$ we take $A(\bos)=([V(\bomega_{i_1,j_1})])$. Assume that we have defined $A(\bos')$ if $\bos'\in\uS\cap \mathbb I_n^{r-1}$. It will be convenient to assume that $[V(\bomega_{i,j})]=0$ if $[i,j]\notin \mathbb I_n$. Suppose that $\bos\in \uS$. For $\bos\in\uS\cap\mathbb I_n^r$, recall that $r_1$ is maximal, so that $\bos(0,r_1)\in\bs^\circ\sqcup\bs$ and define $A(\bos)$ as follows: $$ A(\bos)_{p,\ell}=A(\bos(1,r))_{p-1,\ell-1},\ \ p,\ell >1.$$
In the remaining cases, we set   
\begin{itemize}
\item $A(\bos)_{1,\ell}=A(\bos)_{\ell,1}=0$ if $\bos(0,\ell)$ is not connected or if $\bos(r_1  -1,\ell)\notin\bs^\circ\sqcup\bs$.\end{itemize}
If $\bos(0,\ell)$ is connected and $\ell\le r_1$ then,
\begin{itemize}
\item  $A(\bos)_{1,\ell}=[V(\bomega_{i_1,j_\ell})]$ and $A(\bos)_{\ell,1}=[V(\bomega_{i_\ell,j_1})]$,
\end{itemize}
while if  $\ell>r_1$ and
\begin{itemize}  
\item  $\bos(0,r_1)\in\bs$ with $\bos(r_1-1,\ell)\in\bs^\circ$  then  $A(\bos)_{1,\ell}=0$ and $A(\bos)_{\ell, 1}=[V(\bomega_{i_\ell,j_1})]$,
\vskip6pt
\item  $\bos(0,r_1)\in\bs^\circ$ with  $\bos(r_1-1,\ell)\in\bs$ then $A(\bos)_{1, \ell}=[V(\bomega_{i_1,j_\ell})] $ and $A(\bos)_{\ell,1}=0$. 
\end{itemize}

Let \begin{gather}\label{deta}\Sigma(\bos)=\begin{cases}\{\sigma\in\Sigma_r: a_{\sigma(1),1}\cdots a_{\sigma(r), r}\ne 0\}, &
\bos(0,r_1)\in\bs^\circ,\\
\{\sigma\in\Sigma_r: a_{1, \sigma(1)}\cdots a_{r, \sigma(r)}\ne 0\}, & \bos(0,r_1)\in\bs.
\end{cases}
\end{gather}
Note that \begin{equation}\label{sigma1} \sigma\in\Sigma(\bos)\implies \sigma(1)=p,\ \ 1\le p\le r_1.\end{equation}  

\subsubsection{{\bf Examples}} Suppose that $\bos\in \uS$.\\
\begin{enumerit}
    \item[(i)] If $k=1$ then 
$A(\bos)$ is  the matrix $([V(\bomega_{i_s,j_\ell})])_{1\le s,\ell\le r}$.\\
\item[(ii)] If $\bos\in\mathbb I_n^5\cap\uS$ for some $n\gg 0$  is such that $\bos(0,2)\in\bs^\circ$ and $\bos(p-1,p+2)\notin\bs^\circ\sqcup\bs$ for all $1\leq p\leq 3$, then 
\begin{equation*}\label{blockmatrix}A(\bos) = \begin{bmatrix}
        [V(\bomega_{i_1,j_1})]& [V(\bomega_{i_1,j_2})] &[V(\bomega_{i_1,j_3})] &0 &0\\
     [V(\bomega_{i_2,j_1})]& [V(\bomega_{i_2,j_2})]& [V(\bomega_{i_2,j_3})]&0&0\\   
    0&[V(\bomega_{i_3,j_2})]&[V(\bomega_{i_3,j_3})]&[V(\bomega_{i_3,j_4})]&[V(\bomega_{i_3,j_5})]\\
    0&[V(\bomega_{i_4,j_2})]&[V(\bomega_{i_4,j_3})]&[V(\bomega_{i_4,j_4})]&[V(\bomega_{i_4,j_5})]\\
    0&0&0&[V(\bomega_{i_5,j_4})]&[V(\bomega_{i_5,j_5})]
\end{bmatrix}.\end{equation*}
\item[(iii)] If $\bos\in\mathbb I_n^5\cap\uS$ for some $n\gg 0$ is such that  $\bos(0,2)\in\bs^\circ$, $\bos(1,4)\in\bs$,  $\bos(3,5)\in\bs^\circ$ 
then
\begin{equation*}A(\bos) = \begin{bmatrix}
        [V(\bomega_{i_1,j_1})]& [V(\bomega_{i_1,j_2})] &[V(\bomega_{i_1,j_3})] & [V(\bomega_{i_1,j_4})]&0\\
     [V(\bomega_{i_2,j_1})]& [V(\bomega_{i_2,j_2})]& [V(\bomega_{i_2,j_3})]&[V(\bomega_{i_2,j_4})]&0\\   
    0&[V(\bomega_{i_3,j_2})]&[V(\bomega_{i_3,j_3})]&[V(\bomega_{i_3,j_4})]&0\\
    0&[V(\bomega_{i_4,j_2})]&[V(\bomega_{i_4,j_3})]&[V(\bomega_{i_4,j_4})]&[V(\bomega_{i_4,j_5})]\\
    0&[V(\bomega_{i_5,j_2})]&[V(\bomega_{i_5,j_3})]&[V(\bomega_{i_5,j_4})]&[V(\bomega_{i_5,j_5})]
\end{bmatrix}.\end{equation*}
\vskip12pt
\end{enumerit}
\subsubsection{}\label{stabledef}
We  say that $\bos\in\uS$ is stable if for $1\le p\le r-1$ we have 
 \begin{gather*}i_{p+1}<i_{p-1} \implies ([i_{p+1}, j_{p+1}], [i_{p-1}, j_{p-1}])\in\bs,\\ j_{p-1}<j_{p+1}\implies ([i_{p+1}, j_{p+1}], [i_{p-1}, j_{p-1}])\in\bs^\circ.\end{gather*}
Notice that the conditions obviously hold if $\bos(p-2,p+1)\in\bs^\circ\sqcup\bs$;  otherwise using the definition of $\uS$ we see that $\bos$ is stable if and only if $j_{p+1}<i_{p-1}$ in the first case and $j_{p-1}<i_{p+1}$ in the second case.

Notice that the third example in Section \ref{altsnakesdef} gives an infinite family of stable alternating snakes. \\\\
  For $\sigma\in\Sigma_r$ set $$\sigma(\bos)=\begin{cases}([i_{\sigma(1)},j_1],\cdots, [i_{\sigma(r)}, j_r]),\ \ \bos(0,r_1)\in\bs^\circ\\ ([i_1,j_{\sigma(1)}],\cdots ,[i_r,j_{\sigma(r)}]),\ \ \bos(0,r_1)\in\bs.\end{cases}$$  Our final result on alternating snake modules is the following. In the special case when $\bos\in\bs^\circ$ the result can be deduced from the work of \cite{LM18} by using Schur--Weyl duality and working in large enough rank.
\begin{thm}\label{det}
  Suppose that  $\bos\in\uS$ is stable.\begin{enumerit}
      \item[(i)]
 The following equality holds in $\cal K_0(\mathscr F_n)$:
$$[V(\bomega_{\bos})]=\det A(\bos)=\sum_{\sigma\in\Sigma(\bos)}(-1)^{\sgn(\sigma)} [W(\bomega_{\sigma\bos})].$$
\item[(ii)] If $j_s\ne j_p$ (or $i_s\ne i_p)$   for all $1\le s\neq p\le r$ we have,
$$[V(\bomega_\bos)]=\sum_{\bomegas\in\cal I_n^+} c_{\bomegas,\bomegas_\bos}[W(\bomega)],\ \ c_{\bomegas,\bomegas_\bos}\in\{-1,0,1\}.$$
 \end{enumerit}
\end{thm}
\subsection{Alternating snakes: further examples} \label{exmore}

For $r\ge 1$  set\begin{gather*}P_r=\{(\mu_1,\cdots, \mu_r)\in\mathbb C^r: \mu_s-\mu_{s+1}\in\mathbb Z,\ 1\le s\le r-1\},\\ P_r^\pm=\{(\mu_1,\cdots,\mu_r)\in \mathbb C^r: \mu_s-\mu_{s+1}\ge \mathbb Z_{\geq 0},\ \ 1\le s\le r-1\},
 \\
 P_r^{\rm reg}=\{(\mu_1,\cdots,\mu_r)\in P_r^+: \mu_s\ne \mu_{\ell}\ \ 1\le s,\ell\leq r \},\\
\rho = \left(\frac{r-1}{2}, \frac{r-3}{2},\cdots, \frac{-r+1}{2}\right)\in P_r^{\rm reg}.\end{gather*}
In what follows we will drop the dependence on $r$.
\\\\
For $1\le k\le r$ let   $\bor= (r_0, r_1,\cdots, r_k)\in\mathbb N^{k+1}$ be such that $$r_0=0,\ \  \ r_\ell>r_{\ell-1}+1+\delta_{\ell, 1},\ \ \ 1\le \ell\le k,\ \ \ r_k=r.$$
 We say that  $\mu+\rho=(\mu_1,\cdots,\mu_r)\in P$ is adapted to $\bor$ if, for all appropriate $0\le \ell\le k$, the following hold: 
 \begin{gather}\label{reg1}
 (\mu_{r_{2\ell-1}},\mu_{r_{2\ell}+1}, \mu_{r_{2\ell}+2},\cdots,\mu_{r_{2\ell+1}-2}, \mu_{r_{2\ell+1}-1},\mu_{r_{2\ell+2}})\in P^{\rm reg},\\ \label{reg2} ( \mu_{r_{2\ell+1}}, \mu_{r_{2\ell+1}+1},\cdots,\mu_{r_{2\ell+2}-1}, \mu_{r_{2\ell+2}})\in P^{\rm reg},\\
\label{reg3}
 \mu_1\leq \mu_{r_2-1}, \ \ \mu_{r_{2\ell-1}}\leq \mu_{r_{2\ell+2}-1}, \ \ \mu_{r_{2\ell+1}+1}\leq \mu_{r_{2\ell+4}}.
\end{gather}

\begin{lem}
    Suppose that $\mu+\rho\in P_r\cap\mathbb Z^r$ is adapted to $\bor$. Let $n\in\mathbb N$ and $\lambda+\rho = (\lambda_1,\cdots, \lambda_r)\in  P_r^{\rm reg}\cap\mathbb Z^r$ be such that $$n+1\ge\lambda_s-\mu_\ell\ge  \delta_{\ell,s}, \ \ 1\le \ell,s\le r.$$
Define $\bos\in\mathbb I_n^r$ as follows: 
    \begin{gather*}
        \bos(r_{2\ell}, r_{2\ell+1}-1)=([\mu_{r_{2\ell}+1},\lambda_{r_{2\ell}+1}],\cdots,[\mu_{r_{2\ell+1}-1}, \lambda_{r_{2\ell+1}-1}]),\\
    \bos(r_{2\ell+1}-1 , r_{2\ell+2}) = ([\mu_{r_{2\ell+2}}, \lambda_{r_{2\ell+2}}],\cdots,[\mu_{r_{2\ell+1}}, \lambda_{r_{2\ell+1}}]),
    \end{gather*} for all appropriate $0\le \ell\le k$. 
    Then $\bos$ is a prime stable alternating snake.
\end{lem}
\begin{pf} 
It is clear from our choices that $\bos$ satisfies the first two conditions in the definition of an alternating snake  and, moreover, $\bos$ is connected. To check that part (iii) of Definition \ref{altsnakesdef} holds, notice that for all appropriate $\ell$ and $m$ we have 
\begin{equation}\label{exmu}\bos(r_{2\ell}-1, r_{2\ell+1})\in \bs^\circ, \ \ \ \bos(r_{2\ell+1}-1, r_{2\ell+2})\in \bs, \ \ \ \bos(r_m-2,r_m+1)\notin\bs^\circ\sqcup\bs.\end{equation}
Hence, part (iii) follows by noting that \eqref{reg1}--\eqref{reg3} give
\begin{gather*}
\mu_1\leq \mu_s< \lambda_s< \lambda_{r_1-1}, \ \ s\geq r_1, \ \ s\neq r_2, \\
\mu_{r_{2\ell-1}}\leq \mu_s<\lambda_s\leq \lambda_{r_{2\ell+2}-1}, \ \ s\geq r_{2\ell+1}, \ \ s\neq r_{2\ell+2},\\
    \mu_{r_{2\ell-1}+1}\leq\mu_s <\lambda_s<\lambda_{r_{2\ell}} , \ \ s>r_{2\ell},
\end{gather*}
which also show that $\bos$ is stable. Finally, to prove that $\bos$ is prime, since the $\lambda_p$ are all distinct, \eqref{exmu} implies that it suffices to show that 
$$\mu_{r_{2\ell+1}-1}\neq \mu_{r_{2\ell+2}-1} \ \ {\rm and} \ \ \mu_{r_{2\ell+1}+1}\neq \mu_{r_{2\ell+2}+1}.$$
But this follows by noting that 
$$\mu_{r_{2\ell+1}-1}<\mu_{r_{2\ell-1}}\leq \mu_{r_{2\ell+2}-1} \ \ {\rm and} \ \ \mu_{r_{2\ell+2}+1}> \mu_{r_{2\ell+4}}\geq \mu_{r_{2\ell+1}+1},$$
where we have used \eqref{reg1} for the first inequalities and \eqref{reg3} for the second ones. 
\end{pf}

\subsection{An application to category $\cal O(\lie{gl}_r)$} Let $\lie{gl}_r$ be the Lie algebra of $r\times r$--matrices and let $\lie h$ be the set of diagonal matrices.  We identify $\lie h^*$ with $\mathbb C^r$.  Let $\{\alpha_1,\cdots,\alpha_{r-1}\}\subset P$ be a set of simple roots and  $R^+\subset P$ be the corresponding  set of positive roots for the pair $(\lie{gl}_r, \lie h)$. Fix also a set of coroots $\{h_\alpha:\alpha\in R^+\}\subset \lie h$. \\ 

 Let $\cal O$ be the BGG--category associated to $\lie{gl}_r$ 
In this section we use the Arakawa--Suzuki functor \cite{AS98}, the results of \cite{CPHecke} (see also \cite{GKV}) and Theorem \ref{det} to compute the decomposition of certain (usually not finite--dimensional) irreducible modules in $\cal O$ in terms of Verma modules.
 \subsubsection{} 
 Let $\lie n^+$ be the subalgebra of strictly upper triangular matrices.  The BGG category  $\cal O$ has as objects  finitely generated $\lie g$--modules which are $\lie h$ semi--simple and $\lie n^+$--finite. Among the important objects in $\cal O$ are the Verma module $M(\nu)$ and its irreducible quotient $V(\nu)$ where
 $\nu=(\nu_1,\cdots, \nu_r)\in\lie h^*$.\\\\
 Let $\cal K(\cal O)$ be the Grothendieck group of $\cal O$; it is a free abelian group with basis $[V(\nu)]$, $\nu\in\lie h^*$. The  modules $[M(\nu)]$ are also a basis for  $\cal K(\cal O)$ and hence we can write
 $$[V(\mu)]=\sum_{\nu\in\lie h^*} c_{\mu,\nu} [M(\nu)].$$ It is known that $$c_{\mu,\nu}\ne 0\implies \nu+\rho=w(\mu+\rho),\ \ {\rm for \ some}\ \ w\in \Sigma_r.$$
\subsubsection{The Arakawa--Suzuki functor}\label{as}  We recall some properties of this functor defined in \cite{AS98} and limit ourselves to the case of interest to us. 
\\\\
For $\ell\ge 1$ let $\mathbb H_\ell$ be the degenerate affine Hecke algebra and let ${\rm Rep}(\mathbb H_\ell)$ be the category of finite--dimensional representations. Given $$\lambda+\rho = (\lambda_1,\cdots, \lambda_r)\in P^+\cap \mathbb Z^r, \ \ \  \mu+\rho= (\mu_1,\cdots,\mu_r)\in \mathbb Z^r, \ \ \lambda_i-\mu_i\in \mathbb Z_+, \ \ \ell=\sum_{i=1}^r(\lambda_i-\mu_i),$$ there exists an induced module $M(\lambda,\mu)$ in ${\rm Rep}(\mathbb H_\ell)$ which is called a  standard module. This  module  has a unique irreducible quotient denoted $V(\lambda,\mu)$. 
\\\\
 For $\ell\ge r$ the Arakawa--Suzuki functor $F_\lambda: \cal O\to {\rm Rep}( \mathbb H_\ell)$ is an exact functor satisfying the following:  if $\mu\in P$ is such that  $\lambda_i-\mu_i\in\mathbb Z_+$ for $ 1\le i\le r$ and $\sum_{i-1}^r(\lambda_i-\mu_i)=\ell$ then $$F_\lambda(M(\mu))= M(\lambda,\mu).$$ Otherwise it maps $M(\mu)$ to zero. If in addition we have $\mu(h_\alpha)\le 0$ for all $\alpha\in R^+$ with 
 $ \lambda(h_\alpha)=0$ then
$$ F_\lambda(V(\mu))= V(\lambda,\mu).$$ Otherwise $F_\lambda$ maps $V(\mu)$ to zero.
\subsubsection{From ${\rm Rep}(\mathbb H_\ell)$ to $\mathscr F_n$} \label{cp}It was proved in \cite{Luz89} that ${\rm Rep}(\mathbb H_\ell)$  is equivalent to the category ${\rm Rep}(\hat H_\ell)$ of finite--dimensional representations of the affine Hecke algebra. This category also has a notion of standard modules with unique irreducible quotients and the equivalence preserves standard and irreducible modules. So, we continue to denote the standard and irreducible modules in ${\rm Rep}(\hat H_\ell)$ by $M(\lambda,\mu)$ and $V(\lambda,\mu)$, respectively.\\\\
It was shown in \cite{CPHecke} that there is a functor $F_{\ell,n}:{\rm Rep}(\hat H_\ell)\to \tilde{\mathscr  F_n}$ where $\tilde {\mathscr F_n}$ is the category of finite--dimensional representations of the quantum affine algebra. The functor maps to  the  full subcategory of $\tilde {\mathscr F_n}$ consisting of modules which are subquotients of $\mathbb C_n^{\otimes \ell}$ when regarded as $\bu_q(\lie {sl}_{n+1})$--modules. Moreover it is an equivalence of categories if $\ell\le n$. \\\\
Suppose that  $\ell_1+\ell_2=\ell$; then we have a canonical inclusion of algebras $\hat H_{\ell_1}\times\hat H_{\ell_2}\to\hat H_{\ell}$. Hence if  $M_1, M_2$ are  objects of ${\rm Rep}(\hat H_{\ell_1})$ and ${\rm Rep}(\hat H_{\ell_2})$, respectively we have the corresponding induced module say $M$ for $\hat H_\ell$. The following results were also established in \cite{CPHecke}
\begin{gather*} F_{\ell,n}(M)= F_{\ell_1,n}(M_1)\otimes F_{\ell_2,n}(M_2),\\  
\ell\le n\implies F_{\ell, n}(V(\lambda,\mu))= V(\bomega_{\mu_1,\lambda_1}\cdots\bomega_{\mu_r,\lambda_r}).\end{gather*}
Since $M(\lambda,\mu)$ is the induced module corresponding to one--dimensional representations of $H_{\lambda_s-\mu_s}$, $1\le s\le r$, it follows from the discussion that $$\ell\le n\implies F_{\ell, n}(M(\lambda,\mu))= V(\bomega_{\mu_1,\lambda_1})\otimes\cdots\otimes V(\bomega_{\mu_r,\lambda_r}).$$ We  remind the reader that in the Grothendieck ring  the right hand side has the same
  equivalence class as the corresponding Weyl module.  

  \subsubsection{} We give an application of Theorem \ref{det}.   Suppose that  $\bos\in\uS$ is 
 stable. Choose $\sigma_\bos\in\Sigma_r$ is  such that $$\lambda+\rho=(j_{\sigma_\bos(1)},\cdots,j_{\sigma_\bos(r)})\in P^+,\ \ j_{\sigma_\bos(s)}=j_{\sigma_\bos(p)}, \ \ s<p \implies i_{\sigma_\bos(s)}<i_{\sigma_\bos(p)},$$ and let $\mu+\rho=(i_{\sigma_\bos(1)},\cdots , i_{\sigma_\bos(r)})$. Assume also that $n\gg0$ i.e.,  $$n+1\ge j_{\sigma_\bos(1)}-\min\{i_p: 1\le p\le r\}\ge j_{\sigma_\bos(r)}-\max\{i_p: 1\le p\le r\}\ge 0.$$
 
The following is an immediate consequence of the discussion so far
  and Theorem \ref{det}.
  \begin{prop}\label{applo}
  Retain the notation of this section and     let $ \ell=\sum_{s=1}^r(\lambda_s-\mu_s)$. 
      \begin{enumerit} \item[(i)] We have  $F_{\ell, n}F_\lambda(V(\mu))= V(\bomega_\bos)$.
 \item[(ii)] 
 If $c_{\mu,\nu}\ne 0$ for some  $\nu+\rho=(\nu_1,\cdots, \nu_r)\in P$ then 
 $$c_{\mu,\nu}= \sum_{\sigma\in\Sigma(\bos)}(-1)^{\sgn\sigma}\delta_{\bomegas_{\nu, \lambda},\bomegas_{\sigma(\bos)}},\ \ \bomega_{\nu,\lambda}=\bomega_{\nu_1, j_{\sigma_\bos(1)}}\cdots\bomega_{\nu_r,j_{\sigma_\bos(r)}}.$$ 
If in addition we have $j_s\ne j_p$ for all $1\le s\ne p\le r$ then $c_{\mu,\nu}\in\{-1,0,1\}$ for all $\nu\in P$.
\end{enumerit}\hfill\qedsymbol
\end{prop}
\begin{rem}
    In particular the proposition applies to the pairs $(\lambda+\rho, \mu+\rho)$ defined in Section \ref{exmore}.
\end{rem}

\section{A preliminary collection of  results on alternating snakes and the category $\mathscr F_n$ }\label{known}
In this section we collect together some crucial results on the structure of  $\uS$ and a number of known results on the category $\mathscr F_n$.\\\\
{\em We remind the reader that the element $\bos(p,p')$ was defined in \eqref{spl}, the definition of an element $\bos$ being connected,  prime,  its prime  factors and of being contained in a prime factor was given  in Section \ref{pfsection} and  the  definition of stable in Section \ref{stabledef}. }
 
\subsection{The elements $\Omega(\bos)$ and $\bos^\circ$} \label{indef} 
Given $\bos=([i_1,j_1],\cdots, [i_r,j_r])\in\mathbb I^r_n$ set \begin{gather*}\label{omegaint*}\Omega(\bos)= ([-j_1,-i_1],\cdots, [-j_r,-i_r]),\ \ \ \ \bos^{\circ}=([i_r,j_r],\cdots, [i_1,j_1]).\end{gather*} 
 Clearly $\bomega_{\bos^\circ}=\bomega_\bos$.
  The following is elementary.

\begin{lem}\label{elemalt} Let $\bos\in\uS\cap \mathbb I^r$. 
   \begin{enumerit}
       \item[(i)]  If $p<p'$ then $\bos(p,p')\in \uS$ and, $\bos$ is connected (resp. prime, stable) if and only if $\bos(p,p')$ is connected (resp. prime, stable) for all $0\le p<p'\le r$. Further, for $1\le p\le r$ we have a block decomposition
$$A(\bos)=\begin{bmatrix} A(\bos(0,p))&B_p(\bos)\\ C_p(\bos)& A(\bos(p,r))\end{bmatrix}.$$\\
\item[(ii)]
We have  $\Omega(\bos)\in\uS$ and, for $0\leq \ell_1<\ell_2\leq r$, we have  $\Omega(\bos)(\ell_1,\ell_2)\in\bs^\circ$ if and only if $\bos(\ell_1,\ell_2)\in \bs$. 
Moreover $\Omega(\bos)$ is connected (resp. stable) if and only if $\bos$ is connected (resp.  stable). The prime factors 
 of $\Omega(\bos)$ are obtained by applying $\Omega$ to the prime factors of $\bos$. Further, \begin{gather}\label{Aomega1}A(\Omega(\bos))_{m,\ell}=0\iff A(\bos)_{\ell,m}=0,\\ \label{Aomega2}A(\Omega(\bos))_{m,\ell}=[V(\bomega_{-j_m,-i_\ell})]\iff A(\bos)_{\ell,m}=[V(\bomega_{i_\ell, j_m})].\end{gather} 
 \item[(iii)] We have $\bos^{\circ}\in\uS$ and, $\bos$ is connected if and only if $\bos^\circ$ is connected.
If $\bos^1\vee\cdots\vee\bos^\ell$ is the prime decomposition of $\bos$ then the prime decomposition of $\bos^\circ$ is $(\bos^\ell)^\circ\vee\cdots\vee(\bos^1)^\circ$.\\

\end{enumerit}
\end{lem}
\subsection{} The following elementary result will be used extensively in the paper.\begin{lem} \label{tauprops}  Suppose that $[i_s,j_s]$, $s=1,2,3$, are elements of $\mathbb I_n$ such that
\begin{itemize}
\item the intervals $[i_2,j_2]$ and $[i_3, j_3]$ overlap,
\item the intervals $[i_1,j_1]$ and $[i_s,j_s]$, $s=2,3$, do not overlap.
\end{itemize}
  Then the intervals $[i_1,j_1]$ and $[i_s,j_p] $ with $ \{s,p\}=\{2,3\}$ do not overlap.
 Moreover 
      if  $\bos=([i_1,j_1],\cdots, [i_r, j_r])\in \bs^\circ$ 
      is such that $[i_s, j_s]$ and $[i_{s+1}, j_{s+1}]$ overlap for all $1\le s\le r-1$,
      and  $[i,j]\in\mathbb I_n$ does not overlap $[i_s, j_s]$ for all $1\le s\le r$, then  $[i,j]$ does not overlap  $[i_r, j_1]$.
  \end{lem}
  \begin{pf} Assume without loss of generality  that $i_3<i_2\le j_3<j_2$. 
  There are five possible positions for $j_1$: \begin{gather*}j_1<i_3,\ \  i_3\le j_1<i_2\le j_3<j_2, \ \ i_3<i_2\le j_1<j_3<j_2,\ \\ i_3<i_2\le j_3\le j_1<j_2,\ \ i_3<i_2\le j_3<j_2\le j_1.\end{gather*} The assumptions that $[i_1,j_1]$ and $[i_s, j_s]$ do not overlap for $s=2,3$ imply that we must have the following positions for $i_1$,
\begin{gather*}i_1<j_1<i_3,\ \  i_3\le i_1<j_1<i_2\le j_3<j_2, \ \ i_3<i_2\le i_1<j_1<j_3<j_2,\ \\  i_3<i_2\leq i_1 < j_3=j_1<j_2,\ \ i_3<i_2\le j_3<i_1<j_1<j_2,\\ i_3<i_2\le j_3<j_2<i_1<j_1,\ \
i_1\le i_3<i_2\le j_3<j_2\le j_1,\ \
i_3<i_2\le j_3<i_1<j_2=j_1.\end{gather*} In all cases an inspection shows that $[i_s,j_p]$ and $[i_1,j_1]$ do not overlap for $\{s,p\}=\{2,3\}.$\\\\
For the second assertion of the Lemma, taking the case $r=2$ we have that  $[i,j]$ does not overlap  $[i_2,j_1]$. Proceeding by induction on $r$, assume that $[i_s,j_1]$ does not overlap $[i,j]$ with $s<r$. By our assumptions on $\bos$ we have  $i_{s+1}<i_s\le j_{s+1}<j_1$ and hence the first part of the lemma applies to the intervals $[i_s,j_1], [i_{s+1},j_{s+1}]$ and $[i,j]$ and gives that $[i_{s+1}, j_1]$ and $[i,j]$ do not overlap, which establishes the inductive step and completes the proof. 
  \end{pf}

\subsection{} \label{weylpod} We turn to the representation theory of quantum affine $\lie{sl}_{n+1}$.   Given $\bomega\in\cal  I^+_n$, the Weyl module $W(\bomega)$ is a universal finite--dimensional  cyclic $\widehat\bu_n$--module generated by an $\ell$--highest weight vector $v_{\bomegas}$; this means that, for each $1\leq i\leq n$ and $k\in \mathbb Z$ we have $x_{i,k}^+v_\bomegas=0$ and  $\phi_{i,k}^\pm$ acts on $v_\bomegas$ by a scalar determined by $\bomega$. 
Any  quotient of $W(\bomega)$ is called an $\ell$--highest weight module with $\ell$--highest weight $\bomega$ and it has  a unique irreducible quotient which is isomorphic to $V(\bomega)$.\\\\
For $\bomega,\bomega'\in\cal I_n^+$  the module $V(\bomega\bomega')$ is a subquotient of $V(\bomega)\otimes V(\bomega')$. 
If  $V(\bomega)\otimes V(\bomega')$ and $V(\bomega')\otimes V(\bomega)$ are both quotients of $W(\bomega\bomega')$ then $$V(\bomega)\otimes V(\bomega')\cong V(\bomega\bomega')\cong V(\bomega')\otimes V(\bomega).$$  
The following result was  established in \cite{Ch01} (see also \cite{VV02}) and will play an important role in this paper.
\begin{prop} \label{weylpermute} Suppose that $\bos=([i_1,j_1],\cdots, [i_k, j_k])\in\mathbb  I_n^k$. Then $$W(\bomega_{\bos})\cong V(\bomega_{i_1, j_1})\otimes \cdots\otimes V(\bomega_{i_k, j_k})$$ provided that  for all $1\le p<s\le k$ with $([i_p,j_p],[i_s,j_s])$ connected  we have $i_p+j_{p}\ge i_s+j_s$. In particular,  for   $\bomega,\bomega'\in\cal I_{n}^+ $ we have \begin{equation}\label{weylgroth} [W(\bomega\bomega')]=[W(\bomega)][W(\bomega')],\ \ {\rm and \ so} \ \ \wt_\ell W(\bomega\bomega')=\wt_\ell W(\bomega)\wt_\ell W(\bomega').\end{equation} 
 If $([i_p,j_p],[i_s,j_s])$ are not connected  for all $1\le s,p\le k$ 
 then $$W(\bomega)\cong V(\bomega)\cong V(\bomega_{i_{\sigma(1)}, j_{\sigma(1)}})\otimes \cdots\otimes V(\bomega_{i_{\sigma(k)}, j_{\sigma(k)}}),\ \ \sigma\in\Sigma_k.$$ 
 \end{prop}
 The following is immediate.
\begin{cor}  Suppose that \begin{gather*} \bos'=([i_1',j_1'],\cdots,[i_\ell', j_\ell'])\in\mathbb I_n^\ell,\ \ \bos''=([i_1'',j_1''],\cdots,[i_r'',j_r''])\in\mathbb I_n^r. \end{gather*} Suppose that for every pair $(p,s)$ with $1\le p\le \ell$ and $1\le s\le r$ 
either $ i'_{p}+j'_{p}\ge i''_{s}+j''_{s}$ or  $([i'_{p}, j'_{p}], [i''_{s}, j''_{s}])$ is  not connected. Then $$W( \bomega_{\bos'}\bomega_{\bos''})\cong W(\bomega_{\bos'})\otimes W(\bomega_{\bos''}).$$
 \end{cor}

\subsection{}\label{overlapred} The following assertions are  well--known (see, for instance, \cite{Ch01}) in terms of the old index set $\bvarpi_{i,a}$. We  reformulate that result in the language of this paper. 
\\\\
Suppose that $([i_1,j_1], [i_2,j_2])\in\bs^\circ$ is connected, i.e.,  $i_2<i_1\le j_2<j_1$ and $j_1-i_2\le n+1$. Then
  \begin{gather}
 \label{lrootdrop}
\bomega_{i_1,j_1}^{-1}\bomega_{i_1,j_2}\bomega_{i_2,j_1}\in\wt_\ell V(\bomega_{i_2,j_2}),\\
  \label{domwttp}
     \wt_\ell ^+(W(\bomega_{i_1,j_1}\bomega_{i_2,j_2}))=\{\bomega_{i_1,j_1}\bomega_{i_2,j_2},\ \bomega_{i_1,j_2}\bomega_{i_2,j_1}\}.\\
     \label{fundtpred}
[V(\bomega_{i_1,j_1})][V(\bomega_{i_2,j_2})]=[V(\bomega_{i_1,j_1} \bomega_{i_2,j_2})]+[V(\bomega_{i_1,j_2}\bomega_{i_2,j_1})],\\ \label{dim1tpfund}
\dim (W(\bomega_{i_1,j_1}\bomega_{i_2,j_2}))_{\bomegas_{i_1,j_2}\bomegas_{i_2,j_1}}=1.
 \end{gather}
 If $([i_1,j_1],[i_2,j_2])$ is  not connected, we have 
\begin{gather} \label{fundtpirred}
     \wt_\ell ^+(W(\bomega_{i_1,j_1}\bomega_{i_2,j_2}))=\{\bomega_{i_1,j_1}\bomega_{i_2,j_2}\},\\
 \label{fundtpirr} [V(\bomega_{i_1,j_1})][V(\bomega_{i_2,j_2})]=[W(\bomega_{i_1, j_1}\bomega_{i_2,j_2})]=[V(\bomega_{i_1,j_1}\bomega_{i_2,j_2})].\end{gather}
 In particular, if $([i_1,j_1],[i_2,j_2])$ is connected it follows that  $$[V(\bomega_{i_1,j_2}\bomega_{i_2,j_1})]=[V(\bomega_{i_1,j_2})][V(\bomega_{i_2,j_1})].$$
Notice that \eqref{fundtpirr} implies that $V(\bomega_{i,j})$ is real for all $[i,j]\in\mathbb I_n$. It is also well--known to be prime and that
$$\bomega\in\wt_\ell V(\bomega_{i,j})\iff \dim V(\bomega_{i,j})_\bomegas=1.$$

   \subsection{$\ell$--roots and a partial order on $\cal I_n^+$} \label{lrootdef}For $[i,j]\in\mathbb I_n$ with $0<j-i<n+1$
   set $$\balpha_{i,j}=\bomega_{i,j}\bomega_{i+1,j+1}(\bomega_{i+1,j}\bomega_{i,j+1})^{-1}.$$ Let $\cal Q_n^+$ be the submonoid (with unit) of $\cal I_n$ generated by the elements $\{\balpha_{i,j}: 0<j-i<n+1\}.$ It is well--known that $\cal Q_n^+$ is free on these generators and that if $\gamma\in\cal Q_n^+\setminus\{\bold 1\}$ then $
\gamma^{-1}\notin\cal I_n^+$.\\\\
Define a partial order $\preccurlyeq$ on $\cal I_n^+$ by $\bomega'\preccurlyeq\bomega $ iff $\bomega'=\bomega\balpha^{-1}$ for some $\balpha\in\cal Q_n^+$.
The elements $\{ \bomega_{i,j}: [i,j]\in \mathbb I_n\}$ are minimal with respect to the partial order $\preccurlyeq$. 
   It is well--known (see \cite[Theorem 3]{FR99} for instance) that for $0<j-i<n+1$ we have 
\begin{equation}\label{lrootfund} \bomega_{i,j}\gamma^{-1}\notin\cal I_n^+,\ \ \gamma\in\cal Q_n^+\setminus\{\bold 1\},\ \ {\rm and} \ \ \bomega\in\wt_\ell V(\bomega_{i,j})\implies \bomega=\bomega_{i,j}\ \ {\rm or}\ \ \bomega\preccurlyeq\bomega_{i,j}\balpha_{i,j}^{-1}.\end{equation}
 We isolate the following trivial observation for later use.
\begin{lem}\label{gammcomp}
    Suppose that $\gamma\in\cal Q_n^+\setminus\{\bold 1\}$ and let $$\gamma=\bomega_{i_1,j_1}^{\epsilon_1}\cdots\bomega_{i_r,j_r}^{\epsilon_s}=\balpha_{p_1,\ell_1}\cdots\balpha_{p_s, \ell_s},\ \ \epsilon_m\in\{-1,1\},\ \  1\le m\le r$$ be  reduced expressions for $\gamma$ in the generators $\{\bomega_{i,j}: 0<j-i<n+1\}$ and $\{\balpha_{i,j}: 0<j-i<n+1\}$ respectively. Then \begin{gather*} 
\epsilon_m=1\implies [i_m,j_m]\in\{[p_k,\ell_k], [p_k+1, \ell_k+1]: 1\le k\le s\},\\
\epsilon_m=-1\implies [i_m,j_m]\in\{[p_k+1,\ell_k], [p_k, \ell_k+1]: 1\le k\le s\}.\end{gather*}
In particular, if  $\bomega\in\cal I_n^+$ is such that $\bomega\gamma^{-1}\in\cal I_n^+$ then there exists $1\le k\le s$ such that either $\bomega\bomega_{p_k,\ell_k}^{-1}\in\cal I_n^+$ or $\bomega\bomega_{p_k+1, \ell_k+1}^{-1}\in\cal I_n^+$.
\end{lem}
\subsection{}   The next lemma will be useful in later sections.
   \begin{lem}\label{gammaroot}
   Suppose that $([i_1,j_1], [i_2,j_2])\in\bs$ is connected. Then$$\bomega_{i_1,j_1}\bomega_{i_2,j_2}=\bomega_{i_1,j_2}\bomega_{i_2, j_1}\prod_{i=i_1}^{i_2-1}\prod_{j=j_1}^{j_2-1}\balpha_{i,j}.$$
   \end{lem}
   \begin{pf}
        
For $s\leq j_1$ an induction on $j_{2}-1-j_1$ (with induction beginning when $j_{2}=j_1+1$ by definition of $\balpha_{s,j}$) shows that 
$$\bbeta_s :=\balpha_{s, j_1}\balpha_{s, j_1+1}\cdots\balpha_{s,j_2-1}=
\bomega_{s,j_1}\bomega_{s+1,j_2}\bomega_{s+1,j_1}^{-1}\bomega_{s,j_2}^{-1}.$$
A further induction on $i_2-1-i_1$ along with the  the fact that $i_2\leq j_1$ gives
$$ \bbeta_{i_1}\cdots\bbeta_{i_2-1}=\bomega_{i_1,j_1}\bomega_{i_2,j_2}\bomega_{i_1,j_2}^{-1}\bomega_{i_2,j_{1}}^{-1}$$ and the lemma follows.\end{pf}

 \subsection{}
  The proof of the following can be found in \cite{CM05}:
\begin{prop} \label{weylkl}Let $\bomega\in\cal I_n^+$. 
\begin{enumerit}
\item[(i)] We have  $\dim W(\bomega)_{\bomegas}=1=\dim V(\bomega)_\bomegas$.
\item[(ii)] If  $\bomega'\in \wt_\ell W(\bomega)$ then $\bomega'\preccurlyeq\bomega$.
 In particular in $\cal K_0(\mathscr F_n)$ we have $$[W(\bomega)]= [V(\bomega)]+\sum_{\bomegas'\prec\bomegas} a_{\bomegas' ,\bomegas}[V(\bomega')],\ \ a_{\bomegas' ,\bomegas}\in\mathbb Z_+,$$
 and $a_{\bomegas',\bomegas}\ne 0$ for finitely many choices of $\bomega'$. \end{enumerit}
 \end{prop}
The following is immediate.
\begin{cor} For $\bomega\in\cal I_n^+$, we have,
$$[V(\bomega)]=[W(\bomega)]+\sum_{\bomegas'\prec\bomegas}c_{\bomegas',\bomegas}[W(\bomega')],\ \ c_{\bomegas',\bomegas}\in\mathbb Z$$ and $c_{\bomegas',\bomegas}\ne 0$ for finitely many choices of $\bomega'$.
\end{cor}

\subsection{} We  give a representation theoretic interpretation of the map $\Omega:\mathbb I_n\to \mathbb I_n$ defined in Section \ref{indef}. Define a  homomorphism of groups $\cal I_n\to\cal I_n$ by extending the assignment $\bomega_{i,j}\to\bomega_{-j,-i}$ and continue to denote the homomorphism  by $\Omega$.  Clearly  
$$\Omega(\cal Q_n^+)=\cal Q_n^+\ \ {\rm and}\ \ \bomega'\prec\bomega\iff \Omega(\bomega')\prec \Omega(\bomega).$$
\begin{lem}\label{Omega} There exists a ring  homomorphism $\tilde\Omega: \cal K_0(\mathscr F_n)\to \cal K_0(\mathscr F_n)$ such that
 $$\tilde \Omega([W(\bomega_{\bos})])= [W(\bomega_{\Omega(\bos)})],\ \ \tilde \Omega([V(\bomega_\bos)])=[ V(\bomega_{\Omega(\bos)})],\ \ \bos\in\mathbb I_n^r,\ \ r\ge 1.$$   \end{lem}
\begin{pf} It is known (see \cite{CP95a}, \cite{CP96}) that there exist homomorphisms $\tau_a:\widehat\bu_n\to\widehat\bu_n$, $a\in\mathbb Z$
and $\bar\Omega:\widehat\bu_n\to\widehat\bu_n$ defined on the generators $x_{i,s}^\pm$ for $1\le i\le n$ and $s\in\mathbb Z$  by $$\tau_a( x_{i,s}^\pm) = q^{as}x_{i,s}^\pm, \ \ \ \ \bar\Omega(x_{i,s}^\pm)= - x_{i,-s}^\mp.$$
Denoting  by $\tau_a(V)$ and $\bar\Omega (V)$ the pull back of an object $V$ of $\mathscr F_n$, it was proved in those papers that 
\begin{gather*}
\tau_{n+1}(\bar\Omega (V(\bomega_{i_1,j_1}\cdots\bomega_{i_r,j_r})))\cong  V(\bomega_{-i_1, -j_1+n+1}\cdots\bomega_{-i_r-j_r+n+1}),\\ \tau_{n+1}(\bar\Omega(V_1\otimes V_2))\cong \tau_{n+1}(\bar \Omega(V_2))\otimes \tau_{n+1}(\bar\Omega(V_1)).\end{gather*}
It was also shown that the dual of $V(\bomega_\bos)$ is given by 
$$V(\bomega_\bos)^*\cong V(\bomega_{j_1-n-1,i_1}\cdots\bomega_{j_r-n-1,i_r}). $$ Moreover, since $ (V_1\otimes V_2)^*\cong V_2^*\otimes V_1^*$ for any pair of objects of $\mathscr F_n$ we have 
\begin{gather*}(\tau_{n+1}(\bar\Omega (V(\bomega_{i_1,j_1}\cdots\bomega_{i_r,j_r}))))^*\cong V(\bomega_{-j_1,-i_1}\cdots\bomega_{-j_r,-i_r}),\\
(\tau_{n+1}(\bar\Omega (V_1\otimes V_2)))^*\cong (\tau_{n+1}(\bar\Omega (V_1)))^*\otimes (\tau_{n+1}(\bar\Omega (V_2)))^*.
\end{gather*}
Hence the assignment $\tilde\Omega([V])=[(\tau_{n+1}(\bar\Omega(V)))^*]$ is an endomorphism of the ring $\cal K_0(\mathscr F_n)$ satisfying $\tilde \Omega([W(\bomega_\bos)])=[W(\bomega_{\Omega(\bos)})] $ and $\tilde\Omega([V(\bomega_\bos)])=[V(\bomega_{\Omega(\bos)})]$.   
\end{pf}

\subsection{}\label{mya}  We  reformulate in the  language of intervals  a very special case of a result established in   \cite{MY12}. 
Given $[i,j]\in\mathbb I_n$ let $\mathbb P_{i,j}$ be the set of all functions $g:[0,n+1]\to\mathbb Z$ satisfying the following conditions:
\begin{gather*}g(0)=2j,\ \ g(r+1)-g(r)\in\{-1,1\},\ \ 0\le r\le n,\ \ g(n+1)=n+1+2i.\end{gather*}
 For $g\in\mathbb P_{i,j}$ we have $g(r)-r\in2\mathbb Z$ and we set
\begin{gather*}
\boc_g^\pm =\left\{\left[\frac12(g(r)-r), \frac12(g(r)+r)\right]: 1\le r \le n,\   g(r-1)=g(r)\pm 1=g(r+1)\right\},\\
\bomega(g)=\prod_{[m,\ell]\in\boc_g^+}\bomega_{m,\ell}\prod_{[m,\ell]\in \boc_g^-}\bomega_{m,\ell}^{-1}\in\cal I_n,\\
\boc^\pm_{i,j}=\bigcup_{g\in \mathbb P_{i,j}}\boc_g^\pm.\end{gather*}  The following assertions are well known (see, for instance, \cite[Lemma 5.10]{MY12}): for $g\in \mathbb P_{i,j}$ we have 
\begin{gather}
 [m,\ell]\in\boc_g^-\implies m+\ell>i+j,\ \ 
\label{uclc}[m,\ell]\in\boc_{i,j}^+\iff [m+1, \ell+1]\in\boc_{i,j}^-.\end{gather}
The following result was proved in \cite{MY12}.
\begin{prop}\label{mysnake}
   For  $\bos=([i_1,j_1],\cdots, [i_r,j_r])\in\bs^\circ$,\  let  $\mathbb P_{\bos}$ be the collection of $r$--tuples $(g_1,\cdots, g_r)$ with $g_s\in\mathbb P_{i_s,j_s}$ for $1\le s\le r$ such that
    $$  g_{s}(m)>g_{s+1}(m),\ \ {\rm for\ all}\ \  1\leq s \leq r-1,\ \  0\leq m\leq n+1.$$ Then,
    $$\wt_\ell V(\bomega_\bos)=\{\bomega(g_1)\cdots\bomega(g_k): (g_1,\cdots, g_k)\in\mathbb P_{\bos}\},\ \ \wt_\ell^+ V(\bomega_\bos)=\{\bomega_\bos\}.$$\hfill\qedsymbol
\end{prop}
\subsection{}

We conclude this section with a consequence of Proposition \ref{mysnake}. 
\begin{lem}\label{ellwttau}
   Suppose that $\bos=([i_1,j_1], [i_2,j_2])\in\mathbb I_n^2$ with $i_1+j_1>i_2+j_2$. Then $[i_1,j_1]\in\boc_{i_2,j_2}^-  $ if and only if $([i_1,j_1], [i_2,j_2]) $ is connected. 
\end{lem}
\begin{pf}
 If $([i_1,j_1],[i_2,j_2])$ is connected it follows from \cite[Section 6.4]{MY12a} that there exists a unique $p\in \mathbb P_{i_2,j_2}$ such that $\boc_p^-=\{[i_1,j_1]\}$ and hence $[i_1,j_1]\in \boc_{i_2,j_2}^-$. \\\\
 For the converse, note that 
given $[i,j]\in\mathbb I_n$ and $g\in \mathbb P_{i,j}$ it is immediate from the definition of $\mathbb P_{i,j}$ that
$$-r\leq g(r)-g(0)\leq r \ \ \ {\rm and} \ \ \ r-n-1\leq g(r)-g(n+1)\leq n+1-r,$$
for $1\leq r\leq n$. In particular 
\begin{equation}\label{lowcbound}\max\{2j-r,\ 2i+r\}\leq g(r)\leq \min\{r+2j,\ 2n+2 +2i -r\}.\end{equation}
Equation \eqref{uclc} shows that the the first inequality is strict if $\frac 12[g(r)-r, g(r)+r]\in\boc_g^-$. Taking $[i,j]=[i_2,j_2]$, $r=j_1-i_1$ in \eqref{lowcbound} and using the fact that $[i_1,j_1]\in \boc_{g}^-$ we have
$$\max\{2j_2 - j_1+i_1,\ j_1-i_1+2i_2\}< g(j_1-i_1)=i_1+j_1,$$
and hence $i_2<i_1$ and $j_2<j_1$. Working with the second inequality in \eqref{lowcbound} we have
$$j_1+i_1= g(j_1-i_1)\leq \min \{j_1-i_1+2j_2,\ 2n+2+2i_2 -j_1+i_1\}$$
and hence $i_1\leq j_2$ and $j_1-i_2\leq n+1$ which completes the proof.
\end{pf}

\section{KKOP invariants} \label{sreal}
Throughout the rest of the paper  we shall use freely (see Section \ref{weylpod}) that for all $\bomega_1,\bomega_2\in\cal I_n^+$ the module $V(\bomega_1\bomega_2)$ is a subquotient of $V(\bomega_1)\otimes V(\bomega_2)$. 
\subsection{} \label{kkop} In \cite{KKOP}, the authors defined  for   $\bomega_1, \bomega_2\in\cal I_{n}^+$ a  non--negative integer $\lie d(V(\bomega_1), V(\bomega_2))$ depending on $n$. 
We summarize certain important properties of $\lie d$ in the following proposition. Part (i) follows from the definition of $\lie d$, (ii) is Corollary 3.17 of \cite {KKOP},  (iii) is Proposition 4.2 of \cite{KKOP}, (iv) is Proposition 4.7 of \cite{KKOP} and finally (v) combines Lemma 2.27 and Lemma
2.28 of \cite{KKOP2}.

\begin{prop} Let $\bomega_1, \bomega_2\in\cal I^+_{n}$ and assume that  $V(\bomega_1)$ is a real $\widehat\bu_n$--module. Then,
\begin{enumerit}
    \item [(i)] $\lie d(V(\bomega_1), V(\bomega_2))=\lie d(V(\bomega_2), V(\bomega_1))$.
    \item[(ii)] $\lie d(V(\bomega_1), V(\bomega_2))=0$ if and only if $V(\bomega_1)\otimes V(\bomega_2)$ is irreducible.
    \item[(iii)] For all $\bomega_3\in\cal I_n^+$ we have $$\lie d(V(\bomega_1), V(\bomega_2\bomega_3))\le \lie d (V(\bomega_1), V(\bomega_2))+\lie d(V(\bomega_1), V(\bomega_3)).$$
    \item[(iv)] The module  $V(\bomega_1)\otimes V(\bomega_2)$ has length two if $\lie d(V(\bomega_1), V(\bomega_2))=1$.
    \item[(v)] Suppose that $V(\bomega_1)$ and $V(\bomega_2)$ are both real modules with $\lie d(V(\bomega_1), V(\bomega_2))\le 1$. Then $V(\bomega_1\bomega_2)$ is real.
    \end{enumerit}
\end{prop}
The following is immediate from a repeated application of  part (iii).
\begin{cor}
    Suppose that $\bomega_s\in\cal I_n^+$ for $1\le s\le p$ and assume that $V(\bomega_1)$ is real. 
    Then $$\lie d(V(\bomega_1), V(\bomega_s))=0,\ \ {\rm for\ all}\ \   \ 2\le s\le p\implies \lie d(V(\bomega_1),V(\bomega_2\cdots\bomega_s))=0.$$ 
    \end{cor}
  
 \subsection{} 
The next proposition was proved in \cite{Naoi24}.
 \begin{prop}\label{naoisnake}
    For  $r\geq 2$ let $\bos=([i_1,j_1],\cdots,[i_{r},j_{r}])\in\bs^\circ\sqcup\bs$. Then $$ \lie d(V(\bomega_{i_1,j_1}), V(\bomega_{\bos(1,r)}))\le 1,$$ with equality holding if and only if  $\bos(0,2)$ is connected.
\end{prop}
\subsection{} \begin{prop}\label{real}
For  $\bos\in\uS$ the module $V(\bomega_\bos)$ is real.
\end{prop} 
\begin{pf} We prove the proposition by induction on $r$ with induction beginning when $r=1$. Assume the result holds for $r-1$  and let $r_1\leq r$ be maximal such that $\bos(0,r_1)\in \bs^\circ\sqcup\bs$.  Taking $\bomega_2=\bomega_{\bos(1,r_1)}$ and $\bomega_3=\bomega_{\bos(r_1, r)}$ in Proposition \ref{kkop}(iii)  we have 
$$\lie d(V(\bomega_{i_1,j_1}), V(\bomega_{\bos(1,r)}))\le \lie d (V(\bomega_{i_1,j_1}), V(\bomega_{\bos(1,r_1)})) +\lie d(V(\bomega_{i_1,j_1}), V(\bomega_{\bos(r_1,r)})).$$
By Definition \ref{altsnakea}(iii) we know that the intervals $[i_1,j_1]$ and $[i_p,j_p]$ do not overlap if $p>r_1$. Hence by \eqref{fundtpirr} and Proposition \ref{kkop}(ii) we have $\lie d(V(\bomega_{i_1,j_1}), V(\bomega_{i_p,j_p}))=0$, for $p>r_1$. Then  Corollary \ref{kkop} gives $\lie d(V(\bomega_{i_1,j_1}), V(\bomega_{\bos(r_1,r)}))=0$. \\\\
Since $([i_1,j_1])\vee \bos(1,r_1)=\bos(0,r_1)\in\bs^\circ\sqcup\bs$, Proposition \ref{naoisnake}
gives  $$\lie d(V(\bomega_{i_1,j_1}), V(\bomega_{\bos(1,r_1)}))\le 1\ \ {\rm and\ 
 so}\ \  \lie d (V(\bomega_{i_1,j_1}), V(\bomega_{\bos(1,r)})\le 1.$$
The inductive hypothesis  applies to  Proposition \ref{kkop}(v)  and  so $V(\bomega_\bos)$ is real. This proves  the inductive step and the proof of the proposition is complete.
\end{pf}

\section{Further results on Weyl modules and Proof of Theorem \ref{pres}(i)}
We establish a number of results on Weyl modules which  are needed to prove the main results. At the end of the section we prove Theorem \ref{pres}(i).

\subsection{} Throughout this section we fix an element $\bos=([i_1,j_1],\cdots, [i_r,j_r]) \in\uS$,  $1\le p< r$ and $\epsilon\in\{0,1\}$ such that   \begin{gather*} \bos(p-1,p+1)\ {\rm is\ contained \ in\ a\ prime \ factor\ of}\  \bos,
\\
{\rm and}\ \ \ i_{p+\epsilon}<i_{p+1-\epsilon}\le j_{p+\epsilon}<j_{p+1-\epsilon}.
\end{gather*}
We remind the reader that the notion of prime factor was defined in Section \ref{pfdef}.
We need the following  technical result for our study.

\begin{lem}\label{primefactorprop3} If $1\le s\le r$ is such that
\begin{gather}\label{bdsb}
i_{p+\epsilon}\le i_s<i_{p+1-\epsilon}\le j_{p+\epsilon}\le j_s<j_{p+1-\epsilon}
\\
\label{bdsa}{\rm (resp.}\ \   i_{p+\epsilon}< i_s\le i_{p+1-\epsilon}\le   j_{p+\epsilon}< j_s\le j_{p+1-\epsilon}\ {\rm)},\end{gather} then 
$s=p+\epsilon$ (resp. $s=p+1-\epsilon$).
\end{lem}
 
\begin{pf}
   We prove \eqref{bdsb} when $\epsilon=1$; the proof when $\epsilon=0$ follows by working with $\bos^{\circ}$. The proof of \eqref{bdsa}  follows by working with $\Omega(\bos)$.\\\\
   Under our assumptions   we have that $\bos(p-1,p+1)\in\bs^\circ$ and that 
   $[i_s,j_s]$ and $[i_p,j_p]$ overlap.  
   If $1\le s<p$ then $\bos(s-1,p)\in\bs$ and hence $\bos(s-1,p+1)\notin\bs\sqcup\bs^\circ$. By Definition \ref{altsnakesdef}(iii), the intervals $[i_s, j_s]$ and $[i_{p+1}, j_{p+1}]$ do  not overlap which forces $a_s=a_{p+1}$ for  some $a\in\{i,j\}$. If $s<p-1$ then we have $$i_{p+1}\le i_s<i_{p-1}<i_p\le j_{p+1}\le j_s<j_{p-1}$$ contradicting the fact that $[i_{p-1}, j_{p-1}]$
and $[i_{p+1}, j_{p+1}]$ do not overlap. Hence $s=p-1$ and we are now in the following situation: $$ \bos(p-2,p)\in\bs,\ \ \ a_{p-1}=a_{p+1},\ \ b_{p+1}<b_{p-1}.$$ Definition \ref{pfdef} shows that $\bos$ must have a prime factor of the form $\bos(\ell-1,p)$ for some $\ell\le p$ which contradicts our assumption that $\bos(p-1,p+1)$ is  contained in a prime factor. Hence we have proved that $s\ge p+1$ and \eqref{bdsb} gives $\bos(p-1, s)\in\bs^\circ$. If $s>p+1$ then we would have $i_s<i_{p+1}$ contradicting \eqref{bdsb}. Hence $s=p+1$ and the proof is complete.\end{pf}

\subsection{} Recall from Section \ref{presthm} that 
\begin{equation*}\tau_p\bos=\bos(0,p-1)\vee([i_{p+1} ,j_p], [i_{p}, j_{p+1}])\vee \bos(p+1,r).
\end{equation*}
It follows from \eqref{domwttp} that $\bomega_{i_{p+1},j_p}\bomega_{i_p,j_{p+1}}\in\wt_\ell W(\bomega_{i_p,j_p}\bomega_{i_{p+1}, j_{p+1}})$. Using equations \eqref{qinv} and \eqref{weylgroth} we have\begin{equation}\label{taupwt} \bomega_{\tau_p\bos}\in\wt_\ell (W(\bomega_{\bos(0,p-1)})\otimes W(\bomega_{\bos(p-1,p+1)})\otimes W(\bomega_{\bos(p+1,r)}))=\wt_\ell W(\bomega_\bos).\end{equation}
 Lemma \ref{gammaroot} gives $\bomega_{\tau_p\bos}=\bomega_{\bos}\gamma_{p,p+1}^{-1}$, where 
\begin{equation}\label{gammadef}  \gamma_{p,p+1}=\bomega_{i_p,j_p}\bomega_{i_{p+1},j_{p+1}}(\bomega_{i_p,j_{p+1}}\bomega_{i_{p+1},j_{p}})^{-1}=\prod_{i=i_{p+\epsilon}}^{i_{p+1-\epsilon}-1}\prod_{j=j_{p+\epsilon}}^{j_{p+1-\epsilon}-1}\balpha_{i,j}.\end{equation}
\begin{prop} For $\gamma\in\cal Q_n^+\setminus\{\bold 1\}$ we have  \begin{equation}\label{gammdom} \gamma \preccurlyeq \gamma_{p,p+1} \  \ {\rm and}\ \  \bomega_\bos\gamma^{-1}\in\wt_\ell^+ W(\bomega_\bos)\iff \gamma=\gamma_{p,p+1}.
  \end{equation}
    
\end{prop}
\begin{pf}
It suffices to prove the forward direction; the converse follows from the discussion preceding the proposition.\\\\
Thus let $\gamma\preccurlyeq\gamma_{p,p+1}$ and observe  (see Section \ref{lrootdef}) that  a reduced expression for  $\gamma$ in terms of the generators of $\cal I_n$  must contain $\bomega_{i,j}$ for some $0<j-i<n+1$.
Since $\bomega_\bos\gamma^{-1}\in\cal I_n^+$ we must have  $[i,j]=[i_s, j_s]$ for some $1\le s\le r$.  Lemma \ref{gammcomp} implies that either $\balpha_{i_s,j_s}$ or $\balpha_{i_s-1, j_s-1}$ must occur in a reduced expression for $\gamma$ in terms of the generators of $\cal Q_n^+$.  
Since $\gamma\preccurlyeq \gamma_{p,p+1}$ the same term must also occur on the right hand side of \eqref{gammadef}. Hence  either  $$ i_{p+\epsilon}\le i_s<i_{p+1-\epsilon}\le  j_{p+\epsilon}\le j_s<j_{p+1-\epsilon}\ \ {\rm or}\ \ i_{p+\epsilon}\le i_s-1<i_{p+1-\epsilon}\le j_{p+\epsilon}\le j_s-1\le j_{p+1-\epsilon}.$$
It is immediate from Lemma \ref{primefactorprop3} (equation \eqref{bdsb} or equation \eqref{bdsa}) that 
$s=p$ or $s=p+1$. 
In particular we  have  proved that one of the following must hold:
$$\bomega_{i_p,j_p}\gamma^{-1}\in\cal I_n^+,\ \ {\rm or}\ \ \bomega_{i_{p+1}, j_{p+1}}\gamma^{-1}\in\cal I_n^+,\ \ {\rm or}\ \ \bomega_{i_p,j_p}\bomega_{i_{p+1}, j_{p+1}}\gamma^{-1}\in\cal I_n^+. $$
Equation \eqref{lrootfund} shows that the first two cases cannot happen and so  $\bomega_{i_p,j_p}\bomega_{i_{p+1}, j_{p+1}}\gamma^{-1}\in\cal I_n^+. 
$ \\\\
Next we prove that 
\begin{equation}\label{dim1tp}\dim W(\bomega_\bos)_{\bomegas_\bos\gamma^{-1}}= \dim W(\bomega_{i_p,j_p}\bomega_{i_{p+1},j_{p+1}})_{\bomegas_{i_p,j_{p}}\bomegas_{i_{p+1}, j_{p+1}}\gamma^{-1}}=1.\end{equation}
 For this, we write
$$\bomega_\bos\gamma^{-1}=\bomega_1\cdots\bomega_r,\ \ \ \bomega_s\in\wt_\ell V(\bomega_{i_s,j_s}), \ \ 1\leq s\leq r.$$
We have already proved that  $\balpha_{i_s,j_s}$ can occur in  a reduced expression for $\gamma$ (in terms of the generators of $\cal Q_n$) only if $s=p,p+1$. Hence \eqref{lrootfund} shows that $\bomega_s=\bomega_{i_s,j_s}$ if $s\ne p,p+1$. Along with Proposition \ref{weylkl}(i) and equation \eqref{qinv} it follows  that \begin{gather*}
0\ne \dim W(\bomega_\bos)_{\bomegas_\bos\gamma^{-1}}=\\
\dim W(\bomega_{\bos(0,p-1)})_{\bomegas_{\bos(0,p-1)}}\dim W(\bomega_{i_p,j_p}\bomega_{i_{p+1},j_{p+1}})_{\bomegas_{i_p,j_p}\bomegas_{i_{p+1}, j_{p+1}}\gamma^{-1}}\dim W(\bomega_{\bos(p+1,r)})_{\bomegas_{\bos(p+1,r)}}\\=\dim W(\bomega_{i_p,j_p}\bomega_{i_{p+1},j_{p+1}})_{\bomegas_{i_p,j_{p}}\bomegas_{i_{p+1}, j_{p+1}}\gamma^{-1}},
\end{gather*}
which proves our claim.\\\\
Equations \eqref{domwttp} and \eqref{dim1tpfund} and Lemma \ref{gammaroot} give
$$\gamma=\gamma_{p,p+1}\ \ {\rm and}\ \ 1= \dim W(\bomega_{i_p,j_p}\bomega_{i_{p+1},j_{p+1}})_{\bomegas_{i_p,j_{p+1}}\bomegas_{i_{p+1}, j_p}}=\dim W(\bomega_\bos)_{\bomegas_\bos\gamma_{p,p+1}^{-1}}.$$  Hence \eqref{dim1tp} and so also  the proposition  are proved.
    
\end{pf}

\begin{cor}\label{gammacond} We have $\dim W(\bomega_{\bos})_{\tau_p\bos}=1$ and $\dim V(\bomega_{\bos})_{\tau_p\bos}=0$.
  
\end{cor}

\begin{pf}

Set  \begin{gather*}M:= W(\bomega_{\bos(0,p-1)})\otimes W(\bomega_{i_p,j_p}\bomega_{i_{p+1}, j_{p+1}})\otimes W(\bomega_{\bos(p+1,r)}) ,\ \ \\ U= W(\bomega_{\bos(0,p-1)})\otimes V(\bomega_{i_p,j_p}\bomega_{i_{p+1}, j_{p+1}})\otimes W(\bomega_{\bos(p+1,r)}) .\end{gather*}
   Noting that  $[M]=[W(\bomega_{\bos})]$ the proposition gives  
  $\dim M_{\bomegas_{\tau_p\bos}}=1$. Further, in the course of the proof of the proposition,  we have also proved that $\dim U_{\bomegas_{\tau_p\bos}}$ is equal to $$\dim W(\bomega_{\bos(0,p-1)})_{\bomegas_{\bos(0,p-1)}}\dim V(\bomega_{i_p,j_p}\bomega_{i_{p+1},j_{p+1}})_{\bomegas_{i_p,j_{p+1}}\bomegas_{i_{p+1}, j_p}}\dim W(\bomega_{\bos(p+1,r)})_{\bomegas_{\bos(p+1,r)}}.$$ Hence  \eqref{fundtpred} and \eqref{dim1tpfund} give $\dim U_{\bomegas_{\tau_p\bos}}=0$. 
Since   $V(\bomega_\bos)$ is a further subquotient of $U$ the corollary follows.
    \end{pf}

    \subsection{} We prove some results on tensor product decompositions  of certain  Weyl modules. In all cases it amounts to checking that the conditions of  Corollary \ref{weylpermute} hold.
   
\subsubsection{}  \begin{lem}\label{weylfactor} For $1\le \ell< r$ we have
\begin{gather*}W(\bomega_\bos)\cong \begin{cases} W(\bomega_{\bos(0,\ell)})\otimes W(\bomega_{\bos(\ell,r)}),\ \ \bos(\ell-1, \ell+1)\in\bs^\circ,\\
W(\bomega_{\bos(\ell,r)})\otimes W(\bomega_{\bos(0,\ell)}),\ \ \bos(\ell-1,\ell+1)\in\bs.\end{cases}\end{gather*}
\end{lem}
\begin{pf} 
Let $s\leq \ell <s'$. The definition of $\uS$ gives 
\begin{itemize}
    \item if $\bos(s-1,s')\notin \bs^\circ\sqcup\bs$ then $[i_s,j_s]$ and $[i_{s'},j_{s'}]$ do not overlap;
    \item if $\bos(s-1,s')\in \bs^\circ$ (resp. $\bos(s-1,s')\in \bs$) then $i_s+j_s>i_{s'}+j_{s'}$ (resp. $i_s+j_s<i_{s'}+j_{s'}$).
    \end{itemize}
    An application of  Corollary \ref{weylpermute} gives the result.
\end{pf}
\subsubsection{} Let $r_1\ge 2$ be maximal so that $\bos(0,r_1)\in\bs\sqcup\bs^\circ$.
\begin{lem}\label{weylpfactor}  If $\bos(0,r_1)\in\bs^\circ$ we have
    \begin{equation}\label{taupincl2}
     W(\bomega_{\tau_p\bos})\cong \begin{cases}
     W(\bomega_{\tau_1\bos(0,2)})\otimes W(\bomega_{\bos(2,r)}), & p=1, \  r_1\ge 3\\
        W(\bomega_{\bos(2,r)})\otimes W(\bomega_{\tau_1\bos(0,2)}), & p=1,\   r_1=2\\
        W(\bomega_{i_1,j_1})\otimes W(\bomega_{\tau_{p-1}\bos(1,r)}), & p>1.
         \end{cases}
         \end{equation}
     If $\bos(0,r_1)\in \bs$  we have  a similar statement which is obtained by interchanging the order of the tensor products on the right hand side. 
\end{lem}

\begin{pf} Suppose that $p=1$ and $r_1\ge 3$. If $3\le s\le r_1$ then $i_s+j_s<\min\{i_1+j_1,i_2+j_2\}$. If $s>r_1$ then the intervals $[i_s, j_s]$, $[i_1,j_1]$, $[i_2,j_2]$ satisfy the hypothesis of Lemma \ref{tauprops} and so $[i_1,j_2]$ and $[i_2,j_1]$ do not overlap the interval $[i_s, j_s]$. The hypothesis of Corollary \ref{weylpermute} holds and the first isomorphism is proved.
\\\\
Suppose that $p=1$ and $r_1=2$. Then the intervals $[i_1,j_1]$ and $[i_s, j_s]$ do not overlap if $s\ge 3$. If $\bos(1,s)\notin\bs$ or if $[i_2,j_2]$ and $[i_s, j_s]$ do not overlap then again Lemma \ref{tauprops} shows that $[i_1,j_2]$ and $[i_2,j_1]$ do not overlap the interval $[i_s, j_s]$. Suppose that $\bos(1,s)\in\bs$ for some $s\ge 3$ and that $[i_2,j_2]$ and $[i_s, j_s]$
overlap.  Since $i_2<\min\{i_1,i_s\}$ and $j_2<\min\{j_1,j_s\}$ we have  that either $i_1\leq i_s\le j_2<j_1$ or $
    i_s\le i_1\le j_2<j_s$. Since $[i_1,j_1]$ and $[i_s,j_s]$ do not overlap, either $$i_2<i_1\le i_s\le j_2<j_s\le j_1\ \ {\rm or}\ \ i_2< i_s\le i_1\le j_2<j_1< j_s. $$ An inspection now shows that for $\epsilon\in\{0,1\}$ either 
$[i_{1+\epsilon}, j_{2-\epsilon}]$ does not overlap $[i_s,j_s] $  or  $i_{1+\epsilon}+j_{2-\epsilon}<i_s+j_s$. The hypothesis of Corollary \ref{weylpermute} again holds and so the second isomorphism follows.\\\\
    The proof when $p>1$ is similar.
If  $s\notin\{p,p+1\}$, then  either 
$i_1+j_1>i_s+j_s $ or $ [i_1,j_1]$ and $[i_s,j_s]$  do  not  overlap.
If $\bos(0,p)\notin \bs^\circ\sqcup\bs$ then  $[i_1,j_1]$, $[i_p,j_p]$ and $[i_{p+1},j_{p+1}]$ satisfy the hypothesis of Lemma \ref{tauprops} and hence  the intervals $[i_1,j_1]$ and $[i_{p+1-\epsilon},j_{p+\epsilon}]$ do not overlap, for $\epsilon \in\{0,1\}$. 
If $\bos(0,p)\in \bs^\circ$ then either $\bos(0,p+1)\in \bs^\circ$ or $[i_1,j_1]$ and $[i_{p+1},j_{p+1}]$ do not overlap. In the first case it is clear that $i_1+j_1> i_{p+\epsilon}+j_{p+1-\epsilon}$. In the second case if $[i_1,j_1]$ and $[i_p,j_p]$ do not overlap then an application of Lemma \ref{tauprops} shows that $[i_1,j_1]$ and $[i_{p+\epsilon},j_{p+1-\epsilon}]$ do not overlap for $\epsilon\in\{0,1\}$. If $[i_1,j_1]$ and $[i_{p+1},j_{p+1}]$ do not overlap and $i_p<i_1\le j_p<j_1$ then we must have either $$i_p<i_{p+1}\le i_1\le j_p<j_1\le j_{p+1}\ \ {\rm or}\ \ i_p<i_1\le i_{p+1}\le j_p<j_{p+1}\le j_1.$$
In both cases it is  clear that the hypothesis of Corollary \ref{weylpermute} holds and 
the third isomorphism  is established.
\end{pf}

\subsection{} \label{pres1pf}
The following proposition proves Theorem \ref{pres}(i). 
\begin{prop}\label{etadef}
    There exists a unique (upto scalars) injective  map $\eta_p: W(\bomega_{\tau_p\bos})\to W(\bomega_\bos)$ of $\widehat\bu_n$--modules.
\end{prop}
\begin{pf} It suffices to prove the existence of the map, the uniqueness is immediate from Corollary \ref{gammacond}. 
     The existence of $\eta_p$ is established by induction on $r$, with Section \ref{overlapred} showing  that induction  begins when $r=2$.   \\\\
    For the inductive step suppose that   $p=1$ and let  $\tilde\eta$  be the canonical inclusion  $W(\bomega_{i_1,j_2}\bomega_{i_2,j_1})\hookrightarrow W(\bomega_{i_1,j_1}\bomega_{i_2,j_2})$. If $\bos(0,2)\in\bs^\circ$ (resp. $\bos(0,2)\in\bs$) then Lemma \ref{weylfactor} and Lemma \ref{weylpfactor} show that we have a non--zero injective map of $\widehat\bu_n$--modules $\eta_1: W(\bomega_{\tau_1\bos})\to W(\bomega_{\bos})$ given as follows:
    \begin{gather*}\eta_1=\tilde\eta\otimes\id,\ \ \bos(0,3)\in\bs^\circ\sqcup\bs,\ \ \ \eta_1=1\otimes\tilde\eta,\ \ \bos(0,3)\notin \bs^\circ\sqcup\bs,\\
    {\rm (resp.}\ \ \eta_1=\tilde\eta\otimes\id,\ \ \bos(0,3)\notin\bs^\circ\sqcup\bs,\ \ \ \eta_1=1\otimes\tilde\eta,\ \ \bos(0,3)\in \bs^\circ\sqcup\bs {\rm )}.
    \end{gather*}
     If $p>1$ then $\bos(p-1,p+1)$ is contained in a prime factor of $\bos(1,r)$. Hence the inductive  hypothesis applies and we have an injective map $\tilde \eta_{p-1}:W(\bomega_{\tau_{p-1}\bos(1,r)})\to W(\bomega_{\bos(1,r)})$.  It follows from  the third isomorphism in \eqref{taupincl2} that  $\eta_p:=\id \otimes \tilde\eta_{p-1}$ defines an injective map $W(\bomega_{\tau_p\bos})\to W(\bomega_{\bos})$.
This proves the inductive step and completes the proof of the proposition.
\end{pf}
 We conclude this section with the following observation.
Recall from Section \ref{presthm} that for $1\le \ell \le r_1$ with $\bos(\ell-1,\ell+1)$ contained in a prime factor of $\bos$ we set $M_\ell(\bos)=\eta_\ell(W(\bomega_{\tau_\ell\bos}))$ and  $M_\ell=0$ otherwise. Assume that $\bos(0,r_1)\in\bs^\circ$ and  let $\iota: W(\bomega_{i_1,j_1})\otimes W(\bomega_{\bos(1,r)})\to W(\bomega_{\bos})$ be the (unique up to scalars) isomorphism of Lemma \ref{weylfactor}.
Setting  $K(\bos)= \sum_{\ell=1}^{r_1-1} M_\ell(\bos)$ we see by our construction of $\eta_\ell$ that
\begin{equation}\label{ker}K(\bos)=M_1(\bos)+ \iota(W(\bomega_{i_1,j_1})\otimes K(\bos(1,r)).\end{equation}

 \section{Proofs of Theorem \ref{primefactor} and Theorem \ref{pres}(ii)}  \label{primedecomp}
We assume throughout that $\bos \in\uS$, with $\bos=([i_1,j_1],\cdots, [i_r,j_r])\in\mathbb I_n^r$. 
\subsection{Proof of Theorem \ref{primefactor}(i)}  
Recall that for $1\le p<r$ with $\bos(p-1,p+1)$ contained in a prime factor of $\bos$ we set\begin{equation*}\tau_p\bos=\bos(0,p-1)\vee([i_{p+1} ,j_p], [i_{p}, j_{p+1}])\vee \bos(p+1,r).
\end{equation*}
\begin{prop}\label{lwtprime} 
     Let $\bomega,\bomega'\in\cal I_n^+\setminus\{\bold 1\}$ be such that $\bomega_\bos=\bomega\bomega'$. Suppose that  $1\le p<r$ is such that $\bos(p-1,p+1)$ is contained in a prime factor of $\bos$ and
$\bomega\bomega_{i_p,j_p}^{-1}$ and $\bomega'\bomega_{i_{p+1}, j_{p+1}}^{-1}$ are elements of $\cal I_n^+$. Then $$ \bomega_{\tau_p\bos}\in\wt_\ell(V(\bomega)\otimes V(\bomega'))\setminus\wt_\ell V(\bomega_{\bos}).$$ In particular the module $V(\bomega_\bos)$ is prime    if $\bos\in\uS^{\pr}$
\end{prop}

 \begin{pf} 
Note that Corollary \ref{gammacond} gives $V(\bomega_\bos)_{\bomegas_{\tau_p\bos}}=0$.   Recalling  from \eqref{gammadef} that $\bomega_{\tau_p\bos}=\bomega_{\bos}\gamma_{p,p+1}^{-1}$  we prove that 
\begin{gather}\label{primeweight} {\rm either}\ \ \bomega'\gamma_{p,p+1}^{-1}\in\wt_\ell V(\bomega')\ \ {\rm or}\ \ \bomega\gamma_{p,p+1}^{-1}\in\wt_\ell V(\bomega),\end{gather}
which clearly proves $\bomega_{\tau_p\bos}\in\wt_\ell (V(\bomega)\otimes V(\bomega'))$. 
We prove \eqref{primeweight}  under the assumption that $\bos(p-1,p+1)\in\bs^\circ$; the case $\bos(p-1,p+1)\in\bs$ is obtained by interchanging the roles of $p$ and $p+1$. 
\\\\
Using \eqref{lrootdrop} we have  $\bomega_{i_{p+1}, j_{p+1}}\gamma_{p,p+1}^{-1}\in\wt_\ell V(\bomega_{i_{p+1}, j_{p+1}})$ and so 
$$\bomega'\gamma_{p,p+1}^{-1}\in\wt_\ell (V(\bomega_{i_{p+1}, j_{p+1}})\otimes V(\bomega'\bomega_{i_{p+1}, j_{p+1}}^{-1})).$$
Suppose that $\bomega'\gamma_{p,p+1}^{-1}\in\wt_\ell V(\tilde\bomega)$ where  $V(\tilde\bomega)$ is a subquotient of $V(\bomega_{i_{p+1}, j_{p+1}})\otimes V(\bomega'\bomega_{i_{p+1}, j_{p+1}}^{-1})$. Then $\tilde\bomega\in \wt_\ell^+ W(\bomega')$ and so there exists $\gamma \in \cal Q^+$, with $\gamma\preccurlyeq \gamma_{p,p+1}$, such that
$\tilde\bomega=\bomega'\gamma^{-1}$.  It follows that $$\bomega_\bos\gamma^{-1}=\bomega\tilde\bomega\in\wt_\ell^+ W(\bomega\bomega')=\wt_\ell^+ W(\bomega_\bos).$$
Proposition \eqref{gammacond}  gives that either $\gamma=\bold 1$ or $\gamma=\gamma_{p,p+1}$. In the latter case we have $$\tilde\bomega=\bomega'\gamma_{p,p+1}^{-1}\in\cal I_n^+,\ \ {\rm i.e.}\ \ \bomega'(\bomega_{i_{p+1}, j_{p+1}}\bomega_{i_p,j_p})^{-1}\bomega_{i_p,j_{p+1}}\bomega_{i_{p+1},j_{p}}\in\cal I_n^+.$$
But this is impossible since  $\bomega\bomega_{i_p,j_p}^{-1}\in\cal I_n^+$  and Definition \ref{altsnakea}(i) then  forces $\bomega'\bomega_{i_p,j_p}^{-1}\notin\cal I_n^+$. Hence $\gamma=\bold 1$ proving that $\tilde\bomega=\bomega'$ and \eqref{primeweight} is proved.
\end{pf}

 \subsection{Proof of Theorem \ref{primefactor}(ii)}  Suppose that \eqref{conn} is not satisfied; i.e. there exists $1\le p<r$  such that $\bos(p-1,p+1)$  is not connected. Let $1\leq p_1\leq p<p_2\leq r$ be such that $\bos(p_1-1,p_2)\in \bs^\circ\sqcup\bs$ with  $p_2-p_1$ is maximal. 
Using \cite[Proposition 3.2]{MY12a} we get
$$ V(\bomega_{\bos(p_1-1, p_2)})\cong V(\bomega_{\bos(p_1-1, p)})\otimes V(\bomega_{\bos(p,p_2)}),\ \ {\rm i.e.},\ \ \lie d (V(\bomega_{\bos(p_1-1, p)}), V(\bomega_{\bos(p,p_2)}))=0.$$ By the definition of alternating snakes we have  $[i_s,j_s]$ and $[i_\ell, j_\ell]$ do not overlap if $s<p_1$ and $\ell\ge p+1$ or if $s\le p$ and $\ell>p_2$ we have $$\lie d (V(\bomega_{\bos(0, p_1-1)}), V(\bomega_{\bos(p,r)}))=0=\lie d (V(\bomega_{\bos(p_1-1, p)}), V(\bomega_{\bos(p_2,r)})).$$
An application of  Proposition \ref{kkop} gives \begin{gather*}\lie d (V(\bomega_{\bos(0, p)}), V(\bomega_{\bos(p,r)}))\le \lie d (V(\bomega_{\bos(0, p_1-1)}), V(\bomega_{\bos(p,r)}))+\lie d (V(\bomega_{\bos(p_1-1, p)}), V(\bomega_{\bos(p,r)}))\\=\le \lie d (V(\bomega_{\bos(p_1-1, p)}), V(\bomega_{\bos(p,p_2)}))+\lie d (V(\bomega_{\bos(p_1-1, p)}), V(\bomega_{\bos(p_2,r)}))=0
\end{gather*} and hence $V(\bomega_{\bos(0, p)})\otimes  V(\bomega_{\bos(p,r)})$ is irreducible as needed.
\subsection{Proof of Theorem \ref{primefactor}(iii)} 
Here we are given that $a_{p-1}=a_{p+1}$ for some $a\in\{i,j\}$ and   $\epsilon\in\{0,1\}$ is chosen so that, if $\{a,b\}=\{i,j\}$ then 
$$\bos(p-2,p)\in\bs^\circ\implies b_{p-1+2\epsilon}< b_{p+1-2\epsilon},\ \ \ \bos(p-2,p)\in\bs\implies b_{p+1-2\epsilon}<b_{p-1+2\epsilon}. $$
It follows that  $\bos(p-2, p+1)\notin\bs^\circ\sqcup\bs$
and so the intervals $[i_s, j_s]$ and $[i_\ell, j_\ell]$ do not overlap if $s\le p-1$ or $\ell>p+1$. \\\\
We claim that if $\epsilon=0$ then the intervals $[i_p,j_p]$ and $[i_\ell, j_\ell]$ also do not overlap if $\ell\ge p+2$. This is immediate if $\bos(p-1,\ell)\notin\bs\sqcup\bs^\circ$.
Otherwise, 
suppose that $\bos(p-1,p+2)\in\bs\sqcup\bs^\circ$. If $i_{p+1}=i_{p-1}$ then one of the following holds:
\begin{gather*}i_{p+2}<i_{p+1}=i_{p-1}<i_p\le j_{p+1}< j_{p-1}\ \ {\rm or} \ \ i_{p-1}=i_{p+1}<  j_{p-1}<j_{p+1}<j_{p+2}.\end{gather*} Since $i_{p+1}\le j_{p+2}<j_{p+1}$  the first set of inequalities forces $[i_{p+2}, j_{p+2}]$ and $[i_{p-1}, j_{p-1}]$ to overlap which is a contradiction. Hence the second set of inequalities hold and since $[i_{p+2}, j_{p+2}]$ and $[i_{p-1}, j_{p-1}]$ do not  overlap we get $j_p<j_{p-1}<i_{p+2}\le i_\ell$ for all $\ell\ge p+2$ with $\bos(p-1,\ell)\in\bs^\circ\sqcup\bs$ and the claim is proved in this case. 
If $j_{p-1}= j_{p+1}$ then one of the following holds:
$$i_{p-1}<i_{p+1}\le j_p<j_{p-1}=j_{p+1}<j_{p+2} \ \ {\rm or} \ \ i_{p+2}<i_{p+1}<i_{p-1}\leq j_{p+1}=j_{p-1}.$$
In the first case, since $i_{p+1}<i_{p+2}\leq j_{p+1}$ it follows that $[i_{p-1}, j_{p-1}]$ and $[i_{p+2}, j_{p+2}]$ overlap which is a contradiction. Hence the second set of inequalities hold and, since $[i_{p+2},j_{p+2}]$ and $[i_{p-1},j_{p-1}]$ do not overlap and $i_{p+1}\leq j_{p+2}$, we are forced to have $j_{p+2}<i_{p+1}<i_p$. In particular we get that $j_{\ell}<i_p$ for all $\ell\geq p+2$ with $\bos(p-1,\ell)\in \bs^\circ\sqcup\bs$, thus completing the proof of the claim. 
\\\\
As a consequence of the discussion we have that $$\lie d(V(\bomega_{i_\ell,j_\ell}), V(\bomega_{i_s, j_s}))=0,\ \ s\le p,\ \ \ell\ge p+2.$$
 Proposition \ref{kkop} and its corollary give\begin{gather*} \lie d (V(\bomega_{\bos(0,p)}),V(\bomega_{\bos(p,r)})\le  \lie d(V(\bomega_{\bos(0,p-2)}, V(\bomega_{\bos(p,r)}))+ \lie d (V(\bomega_{\bos(p-2,p)}),V(\bomega_{\bos(p,r)}))
\\ =\lie d (V(\bomega_{\bos(p-2,p)}),V(\bomega_{\bos(p,r)}))\le \lie d (V(\bomega_{\bos(p-2,p)}),V(\bomega_{\bos(p+1, r)}))+\lie d (V(\bomega_{\bos(p-2,p)}),V(\bomega_{\bos(p,p+1)}))
\\
=\lie d(V(\bomega_{\bos(p-2,p)}),V(\bomega_{i_{p+1},j_{p+1}})).\end{gather*} 

In particular, this reduces the proof of part (iii) to the case when $r=3$; hence we assume from now on that $\bos=([i_1,j_2], [i_2,j_2], [i_3, j_3])$.\\\\
Suppose that $\bos(0,2)\in \bs^\circ$; then we have $a_1=a_3$ and $b_1<b_3$. Note that 
$V(\bomega_{i_3,j_3})\otimes V(\bomega_{\bos(0,2)})$ is $\ell$--highest weight by Lemma \ref{weylfactor}. Hence by \cite[Corollary 3.16]{KKKO15} it suffices to prove that if $V(\bomega)$ is in the socle of this tensor product then $\bomega=\bomega_{\bos}$.   Using 
\cite[Lemma 1.3.4]{BC23} and Proposition \ref{mysnake} we see that 
$$\bomega=\bomega_{i_3,j_3}\bomega(g_1)\bomega(g_2),\ \  (g_1,g_2)\in\mathbb P_{\bos(0,2)}.$$
If $\bomega\ne\bomega_{\bos}$ there exists $m,s$ with $\{m,s\}=\{1,2\}$ satisfying  
$$\boc_{g_s}^-=\{[i_3,j_3]\}, \ \ \ g_s(j_3-i_3)=i_3+j_3, \ \ \ \bomega(g_m)=\bomega_{i_m,j_m}.$$
Since $[i_1,j_1]$ and $[i_3,j_3]$  do not overlap, Lemma \ref{ellwttau} forces $s=2$ and $m=1$.  Proposition \ref{mysnake} gives 
\begin{gather}\label{pathk2} g_2(j_{3}-i_{3})= i_{3}+j_{3} < g_1(j_3-i_3)=j_{1}+i_{1}+|j_{1}-i_{1}-j_{3}+i_{3}|,\end{gather}
or equivalently using the fact that  $b_1\le b_3$, $$a_3+b_3<a_1+b_1+b_3-b_1,\ \ {\rm i.e.} \ \ a_3<a_1$$ contradicting our assumption that $a_1=a_3$.
By Lemma \ref{Omega} we have $$[V(\bomega_{\Omega(\bos)})]=[V(\bomega_{\Omega(\bos(0,2))})]
[V(\bomega_{\Omega(\bos(2,3))})],$$
and hence the  irreducibility of $V(\bomega_{i_3, j_3})\otimes V(\bomega_{\bos(0,2)})$ follows in the case 
when $\bos(0,2)\in\bs$.\\\\
This completes the proof of part (iii) of the theorem when $\epsilon=0$. If $\epsilon=1$ then working $\bos^{\circ}$ gives the result.

\subsection{Proof of Corollary \ref{tensorfactor}} If $r=3$ there is nothing to prove and so we assume that $r\ge 4$. Assume also that $\bos\notin\uS^{\pr}$ and let $1\le p<r$ be such that $\bos(0,p)$ is a prime factor of $\bos$. 
 By parts (ii) and (iii) of Theorem \ref{primefactor} we have that $$V(\bomega_{\bos})\cong V(\bomega{_{\bos(0,p)})\otimes V(\bomega_{\bos(p,r)}}).$$ Since $\bos(p,r)$ is a concatenation of the other prime factors of $\bos$, the first statement of the corollary is now immediate by a straightforward induction on $r$. For the second statement,
suppose that $$V(\bomega_\bos)\cong V(\bomega_1)\otimes V(\bomega'), \ \ \ \bomega', \bomega_1\in\cal I_n^+\setminus\{\bold 1\} ,\ \ \bomega_1\bomega_{i_1,j_1}^{-1}\in\cal I_n^+,$$ and  $V(\bomega_1)$ is prime. Let $1\leq p'\le r$ be maximal so that $\bomega_1\bomega_{i_s,j_s}^{-1}\in\cal I_n^+$ for all $1\le s\le p'$. If $p'<p$ then $\bomega'\bomega_{i_{p'+1}, j_{p'+1}}^{-1}\in\cal I_n^+$. Proposition \ref{lwtprime} applies since $\bos(p'-1,p'+1)$ is contained in the prime factor $\bos(0,p)$ and gives that  $V(\bomega_1)\otimes V(\bomega')$ is reducible contradicting our assumptions. Hence $p'=p$ and  
$\bomega_1=\bomega_{\bos(0,p)}\bomega_1'$.\\\\ Suppose that $\bomega_1'\bomega_{i_{p_1},j_{p_1}}^{-1}\in\cal I_n^+$ for some $p_1\ge p+1$ and $p_1$ is minimal with this property. If $\bos(p_1-2,p_1)$ is contained in a prime factor of $\bos$ then Proposition \ref{lwtprime} again shows  that $V(\bomega_1)\otimes V(\bomega')$ is reducible. Hence there exists $p_2\ge p_1$ such that $\bos(p_1-1,p_2)$ is a prime factor of $\bos$. The same  arguments now  show that $\bomega_1\bomega_{i_m,j_m}^{-1}\in\cal I_n^+$ for all $p_1\le m\le p_2$. Repeating we find that $\bomega_1=\bomega_{\bos(0,p)}\bomega_{\bos(p_1-1,p_2)}\cdots\bomega_{\bos(p_{m-1}-1, p_m)}$ where $\bos(p_{\ell-1}-1, p_\ell)$ for $1\le \ell\le m$ are all prime factors of $\bos$. By part (iii) of Theorem \ref{primefactor} the tensor product of the modules associated to any subset of the prime factors of $\bos$ is irreducible and so we have $$V(\bomega_1)\cong V(\bomega_{\bos(0,p)})\otimes V(\bomega_{\bos(p_1-1,p_2)}\cdots\bomega_{\bos(p_{m-1}-1, p_m)}), $$ contradicting our assumption that $V(\bomega_1)$ is prime. Hence $$\bomega_1=\bomega_{\bos(0,p)},\ \ V(\bomega_{\bos(p,r)})\cong V(\bomega').$$ The second assertion of the  corollary is now immediate by an induction on $r$.

\subsection{Proof of Theorem \ref{pres}(ii)} 
Let $\pi: W(\bomega_\bos)\to V(\bomega_\bos)\to 0 $ be the canonical map of $\widehat\bu_n$--modules.
By Theorem \ref{pres}(i)  there exists a unique (upto scalars) non--zero injective map $\eta_p: W(\bomega_{\tau_p\bos})\to W(\bomega_\bos)$ if $\bos(p-1,p+1)$ is contained in a prime factor of $\bos$. 
Let $M_p(\bos)$ be the image of $\eta_p$ if $\bos(p-1,p+1)$ is contained in a prime factor of $\bos$ and otherwise $M_p(\bos)=0$. Recall also from Section \ref{pres1pf} that we set $$K(\bos)=\sum_{p= 1}^{r-1}M_p(\bomega_\bos).$$
It  follows from Proposition \ref{lwtprime} that $\pi(M_p(\bos))=0$ and hence we have a surjective map
$$ \frac{W(\bomega_{\bos})}{K(\bos)}\to V(\bomega_\bos)\to 0.$$ 
We prove that this map is an isomorphism proceeding by induction on $r$. Section \ref{overlapred} (see \eqref{fundtpred}) shows that induction begins at $r=2$. \\\\
Assume  that $\bos(0,2)\in \bs^\circ$.  If $\bos(0,2)$ is not contained in a  prime factor of $\bos$  then $M_1(\bos)=0$ by definition. By Lemma \ref{weylfactor} and Theorem \ref{primefactor}  we have $$W(\bomega_\bos)\cong V(\bomega_{i_1,j_1})\otimes W(\bomega_{\bos(1,r)}),\ \ V(\bomega_\bos)\cong V(\bomega_{i_1,j_1})\otimes V(\bomega_{\bos(1,r)}). $$ 
By the  inductive hypothesis we have a short exact sequence $$0\to \sum_{p=1}^{r-2}M_p(\bos(1,r))\to W(\bomega_{\bos(1,r)})\to V(\bomega_{\bos(1,r)})\to 0.$$  Tensoring with $V(\bomega_{i_1,j_1})$ on the left and using \eqref{ker} gives the inductive step.\\\\
Assume now that $\bos(0,2)$   is contained in a prime factor of $\bos$. Then  the inductive hypothesis gives  a short exact sequence 
 $$0\to \sum_{p=2}^{r-1} M_p(\bos)\to W(\bomega_{i_1,j_1})  \otimes W(\bomega_{\bos(1,r)})  \to V(\bomega_{i_1,j_1})\otimes V(\bomega_{\bos(1,r)})\to 0.$$
 By Proposition \ref{lwtprime}  the module $V(\bomega_{i_1,j_1})\otimes V(\bomega_{\bos(1,r)})$ is reducible. Let  $1<r_1\leq r$ be  maximal such that $\bos(0,r_1)\in\bs^\circ$. Since $[i_1,j_1]$ and $[i_s,j_s]$ do not overlap if $s>r_1$ it follows from Proposition \ref{kkop} and Proposition \ref{naoisnake} that $$0<\lie d(V(\bomega_{\bos(1,r)}),  V(\bomega_{i_1,j_1}))\le \lie d(V(\bomega_{\bos(1,r_1)}), V(\bomega_{i_1,j_1}))=1.$$
  Hence $V(\bomega_{i_1,j_1})\otimes V(\bomega_{\bos(1,r)})$ has length two by Proposition \ref{kkop}(iv). Proposition \ref{lwtprime} and Theorem \ref{pres}(i) show that the composite map $$\eta_1: W(\bomega_{\tau_1\bos})\to W(\bomega_\bos)\to V(\bomega_{i_1,j_1})\otimes V(\bomega_{\bos(1,r)})$$ is non--zero while the further composite to $V(\bomega_\bos)$ is zero. Hence we have the following,
 \begin{gather*}0\to V(\bomega_{\tau_1\bos})\to V(\bomega_{i_1,j_1})\otimes V(\bomega_{\bos(1,r)})\to V(\bomega_\bos)\to 0,\\\\
  \dfrac{\sum_{p=1}^{r-1}M_p(\bomega_\bos)}{\sum_{p=2}^{r-1}M_p(\bomega_\bos)}\hookrightarrow V(\bomega_{i_1,j_1})\otimes  V(\bomega_{\bos(1,r)})\to \frac{W(\bomega_\bos)}{\sum_{p=1}^{r-1}M_p(\bomega_\bos)}\to 0.\end{gather*} It is immediate that $$V(\bomega_\bos)\cong \frac{W(\bomega_\bos)}{\sum_{p=1}^{r-1}M_p(\bomega_\bos)}.$$ If $\bos(0,2)\in\bs$ the proof is identical if one switches the order of the tensor products.

\section{Proof of Theorem \ref{det}}
The proof of Theorem \ref{det} is fairly involved and it requires  additional representation theory.  This theory is interesting in its own right since (see Proposition \ref{idp}) it involves certain cluster type identities. 
We also need several results on the matrix $A(\bos)$ where $\bos\in\uS$ is stable. In turn these depend on a detailed understanding of the structure of alternating snakes. We begin  by stating certain key results whose proofs are given in subsequent sections.  Assuming these results we complete  the proof of    Theorem \ref{det}. \\\\
Throughout this section we fix an  element $\bos\in\uS$. Writing $\bos=([i_1,j_1], \cdots,[i_r,j_r])\in\mathbb I_n^r$ we let $2\le r_1\le r$ be maximal such that $\bos(0,r_1)\in \bs^\circ\sqcup\bs$. We also set $$\bos_1=\bos(1,r),\ \ \bos_p=\begin{cases}([i_1,j_2],\cdots, [i_{p-1}, j_p])\vee\bos(p,r),\ \ \bos(0,r_1)\in\bs^\circ,\\
([i_2,j_1],\cdots, [i_p,j_{p-1}])\vee\bos(p,r),\ \ \bos(0,r_1)\in\bs,
\end{cases}\ \ 2\le p\le r_1.$$
\\
{\em It is convenient to adopt the convention that $$[V(\bomega_{i,j}\bomega)]=0,\ \ {\rm for\ all}\ \  \bomega\in\cal I_n^+,\ \  {\rm if}\ \ j-i<0\ \ {\rm or}\ \  j-i>n+1.$$}

\subsection{}\label{spdef} 
Our first result establishes an identity  in $\cal K_0(\mathscr F_n)$. 
Recall that  $\bos\in\uS$ is connected if and only if $0\le j_{s+1}-i_s,\ j_s-i_{s+1}\le n+1$, for $1\leq s<r$. 
\begin{prop}\label{altsumid}
    Assume that  $\bos\in\uS^{\pr}$ and it is stable. Then   the following equality in $\cal K_0(\mathscr F_n)$:
  $$[V(\bomega_\bos)]=\sum_{p=1}^{r_1}(-1)^{p+1}[V(\bomega_{\bos_p})]\begin{cases}[V(\bomega_{i_p,j_1})],\ \ \bos(0,r_1)\in\bs^\circ,\\
[V(\bomega_{i_1,j_p})],\ \ \bos(0,r_1)\in\bs.\end{cases}$$
\end{prop}
\subsection{} Our next result studies the elements $\bos_p$, $1\le p\le r_1.
$
\begin{prop}\label{sps}
  Suppose that $\bos\in\uS^{\pr}$ is stable. Then  $\bos_p$   is a stable element of $\uS$  for all $1\le p\le r_1$. 
\end{prop}
\begin{rem}
    In view of Proposition \ref{sps} we have that the matrix $A(\bos_p)$ is defined.
\end{rem}

\subsection{} If $\bos(0,r_1)\in \bs^\circ$ (resp. $\bos(0,r_1)\in \bs$) let $A_p(\bos)$, $1\le p\le r_1$, be the matrix  obtained from $A(\bos)$ by dropping the first column (resp. first row)  and the $p$--th row (resp. $p$--the column).
\begin{prop} Assume  that $\bos\in\uS$ is stable. \label{matrixprop}
\begin{enumerit}
\item[(i)] If   $\bos\in\uS^{\pr}$ then $$\det A(\bos_p)=\det A_p(\bos),\ \ 1\le p\le r_1.$$
    \item[(ii)] Suppose that $\bos\in\uS\setminus\uS^{\pr}$ and  that $\bos(0,\ell)$ is a prime factor of $\bos$ for some $1\le \ell< r$. Then,$$\det A(\bos)=\det A(\bos(0,\ell))\det A(\bos(\ell,r)).$$
    
\end{enumerit}

\end{prop}

\subsection{Proof of Theorem \ref{det}(i)} By Proposition \ref{weylpermute} we have that $$[W(\bomega)]=[V(\bomega_{m_1,\ell_1})]\cdots[V(\bomega_{m_s,\ell_s})]\ \ {\rm if}\ \ \bomega=\bomega_{m_1,\ell_1}\cdots\bomega_{m_s,\ell_s}.$$ Hence \eqref{deta} gives  \begin{equation}\label{altsum} \det A(\bos)=\sum_{w\in\Sigma(\bos)}(-1)^{\sgn(w)} [W(\bomega_{w\bos})].\end{equation}
We prove by induction on $r$ that $$[V(\bomega_\bos)]=\det A(\bos).$$
Induction clearly begins at $r=1$ and we assume  that  the result holds for $r-1$.
We  prove  the inductive step when $\bos(0,r_1)\in\bs^\circ$.
The case when $\bos(0,r_1)\in\bs$ follows since an  application of  Lemma \ref{elemalt}(ii) and  Lemma \ref{Omega} 
gives $$\det(A(\Omega(\bos)))=[\tilde\Omega(V(\bomega_{\bos}))]=[V(\bomega_{\Omega(\bos)})].$$  The assumption that $\bos(0,r_1)\in\bs^\circ$   gives $$A (\bos)_{s,1}=0,\ \ s\ge r_1+1,\ \ A(\bos)_{s,1}=[V(\bomega_{i_s,j_1})],\ \ 1\le s\le r_1. $$   Hence  
\begin{equation}\label{detr1}\det A(\bos)=\sum_{p=1}^{r_1}(-1)^{p+1}[V(\bomega_{i_p,j_1})]\det A_p(\bos).\end{equation}
If $\bos\in \uS^{\pr}$, by Proposition \ref{sps} we have that $\bos_p\in\uS$ is stable and hence the inductive hypothesis gives  $[V(\bomega_{\bos_p})]=\det A(\bos_p)$. Then Proposition \ref{altsumid}, Proposition \ref{matrixprop}(i) and equation \eqref{detr1} give
 $$
[V(\bomega_{\bos})]=\sum_{p=1}^{r_1}(-1)^{p+1} [V(\bomega_{i_p,j_1})]\det A(\bos_p)= \det A(\bos).
$$ If $\bos\notin\uS^{\pr}$ then choose $\ell<r$ such that $\bos(0,\ell)$ is a prime factor of $\bos$. The inductive hypothesis applies to $[V(\bomega_{\bos(0,\ell )})]$ and $[V(\bomega_{\bos(\ell,r)})]$. Theorem \ref{primefactor} and Proposition \ref{matrixprop}(ii) give
$$[V(\bomega_\bos)]=[V(\bomega_{\bos(0,\ell )})][V(\bomega_{\bos(\ell,r)})]=\det A(\bos(0,\ell))\det A(\bos(\ell,r))= \det A(\bos)$$ and the inductive step is established.

\subsection{Proof of Theorem \ref{det}(ii)} 
It follows from Corollary \ref{weylkl} and Theorem \ref{det}(i) that we  can write \begin{equation}\label{sigmachoice}[V(\bomega_\bos)]=\sum_{\bomegas\in\cal I_n^+}c_{\bomegas,\bomegas_{\bos}}[W(\bomega)],\ \ \ \ c_{\bomegas,\bomegas_\bos}=\sum_{\sigma\in\Sigma(\bos)}(-1)^{\sgn\sigma}\delta_{\bomegas,\bomegas_{\sigma(\bos)}}.\end{equation}
It is convenient to adopt the convention  that $c_{\bomegas,\bomegas'}=0$ if $\bomega$ or $\bomega'$ are not in $\cal I_n^+$.  Define $$\supp\bomega_\bos=\{\bomega\in\cal I_n^+: c_{\bomegas,\bomegas_\bos}\ne 0 \}.$$ Then \eqref{sigmachoice} shows that $\bomega\in\supp\bomega_\bos$ only if $\bomega=\bomega_{i_1, j_{\sigma(1)}}\cdots\bomega_{i_r, j_{\sigma(r)}}$ for some $\sigma\in\Sigma_r$. \\\\ 
 We prove that $c_{\bomegas,\bomegas_\bos}\in\{-1,0,1\}$ if $j_s\ne j_\ell$ for $1\le s\ne \ell\le r$ by induction on $r$ with induction beginning when $r=1$. The case when $i_s\ne i_\ell$ for all  $1\le s\ne \ell \le r$ follows by working with  $\Omega(\bos)$. 
\\\\ 
Suppose that $\bos\in\uS\setminus \uS^{\pr}$ and let $\bos(0,\ell)$ be  a prime factor of $\bos$ for some $1\le \ell< r$. By Theorem \ref{primefactor} we have 
\begin{gather*}[V(\bomega_{\bos})]=[V(\bomega_{\bos(0,\ell)})][V(\bomega_{\bos(\ell,r)})],\ \ 
 {\rm and \ 
  so}\ \ c_{\bomegas,\bomegas_\bos}=\sum_{\bomegas_1\in\cal I_n^+}c_{\bomegas_1,\bomegas_{\bos(0,\ell)}}
c_{\bomegas\bomegas_1^{-1},\bomegas_{\bos(\ell,r)}}.\end{gather*}
Since $j_s\ne j_p$ for $1\le s\neq p\le r$, it  is  clear that if $\{\bomega_1,\bomega_2\}\subset \supp\bomega_{\bos(0,\ell)}$ and $\{\bomega_1', \bomega_2'\}\subset \supp\bomega_{\bos(\ell, r)}$ are such that $\bomega_1\bomega_1'=\bomega_2\bomega_2'$ then $\bomega_1=\bomega_2$ and $ \bomega_1'=\bomega_2'$. In other words $c_{\bomegas_1,\bomegas_{\bos(0,\ell)}}
c_{\bomegas\bomegas_1^{-1},\bomegas_{(\bos(\ell,r)}}\ne 0$ for at most one choice of $\bomega_1$ and the inductive step follows.
\\\\ 
It remains to   prove the inductive step when $\bos \in \uS^{\pr}$. 
 Proposition \ref{altsumid} gives,
\begin{equation}\label{cs} c_{\bomegas,\bomegas_\bos}=\begin{cases}\sum_{p=1}^{r_1}(-1)^{p+1}c_{\bomegas\bomegas_{i_p,j_1}^{-1}, \bomegas_{\bos_p}},\ \ \bos\in\bs^\circ,\\
\sum_{p=1}^{r_1}(-1)^{p+1}c_{\bomegas\bomegas_{i_1,j_p}^{-1}, \bomegas_{\bos_p}},\ \ \bos\in\bs.
\end{cases}
\end{equation}
 Suppose that  $\bos(0,r_1)\in \bs^\circ$. If $\bomega\in\supp\bomega_\bos$ 
 then by \eqref{sigmachoice} we can choose  $\sigma\in\Sigma(\bos)$ with $\bomega=\bomega_{i_{\sigma(1)},j_1}\cdots\bomega_{i_{\sigma(r)},j_r}.$ Recall from \eqref{sigma1}  that $1\le \sigma(1)\le r_1$.
 In particular since $j_s\ne j_\ell$ if $1\le s\ne \ell \le r$ this means that $\bomega\bomega_{i_p,j_1}^{-1}\in\cal I_n^+$  if and only if  $p=\sigma(1)$. By Proposition \ref{sps}, the induction hypothesis applies to $\bos_p$ and gives
 $$c_{\bomegas, \bomegas_\bos}
=(-1)^{\sigma(1)+1}c_{\bomegas\bomegas_{i_{\sigma(1)}, j_1}^{-1},\bomegas_{\bos_{\sigma(1)}}}\in\{-1,1\}.$$
If $\bos(0,r_1)\in\bs$ the proof is slightly different since we are not assuming that $i_s\ne i_\ell$ if $1\le s\ne \ell\ne r$.
Let $\bomega'\in\supp\bomega_{\bos_p}$ for some  $1\le p\le r_1$ and regard $\Sigma(\bos_p)$ as the set of permutations of $\{1,2,\cdots, r\}\setminus \{p\}$. If $r_1>2$ then $\bos_p(0,r_1-1)\in\bs$ and by \eqref{sigmachoice} we can choose  $\sigma\in\Sigma(\bos_p)$ such that $$\bomega'=\bomega_{\sigma\bos_p}=\bomega_{i_2, j_{\sigma(1)}}\cdots \bomega_{i_r, j_{\sigma(r)}},\ \ \  \sigma(s)\in\{1,2,\cdots, r\}\setminus\{ p\}.$$
Since $j_s\ne j_\ell$ for all $1\le s\ne\ell\le r$  we have $$c_{\bomegas',\bomegas_{\bos_p}}\ne 0\implies c_{\bomegas',\bomegas_{\bos_\ell}}=0,\ \ 1\le p\ne \ell\le r_1.$$
Hence if $\bomega\in\supp\bomega_\bos$,   there exists a unique $p$ such that $\bomega\bomega_{i_1,j_p}^{-1}\in\supp \bomega_{\bos_p}$ and the inductive hypothesis applied to $\bos_p$  gives $c_{\bomegas,\bomegas_{\bos}}\in\{-1,1\}$. If $r_1=2$ then $\bos_p\in\bs^\circ$ for $p=1,2$ and so if $c_{\bomegas',\bomegas_{\bos_1}}\ne 0$ (resp. $c_{\bomegas',\bomegas_{\bos_2}}\ne 0$) there exists $\sigma\in\Sigma(\bos_1)$ (resp. $\sigma\in\Sigma(\bos_2)$)  such that $$\bomega'=\bomega_{i_{\sigma(2)},j_2}\cdots\bomega_{i_{\sigma(r)}, j_r},\ \ {\rm ( resp.} \ \bomega'=\bomega_{i_{\sigma(2)},j_1}\bomega_{i_{\sigma(3), j_3}}\cdots\bomega_{i_{\sigma(r)}, j_r}\ {\rm )}. $$
It again follows that at most one  of $c_{\bomegas',\bomegas_{\bos_p}}\ne 0$ for $p=1,2$ and the inductive step is complete if $\bos$ is prime. The proof of the theorem is complete.

\section{Structure of alternating snakes and Proof of Proposition \ref{sps}}
For our further study, we need several results on the structure of alternating snakes. We collect all of them in this section. We warn the reader that the  proofs are tedious and  the remaining sections of the paper can be read independent of the proofs given here. Throughout this section we fix $\bos= ([i_1,j_1],\cdots, [i_r,j_r])\in \uS$ with $r\geq 2$ and let $2\leq r_1\leq r$ be maximal such that $\bos(0,r_1)\in \bs^\circ\sqcup\bs$.
 \subsection{}  We study stable elements of $\uS$.
  \begin{prop}\label{structure}
    Suppose that $\bos\in \uS$ is  stable and connected and that  $2\leq r_1<r$.   
    Then  $i_{r_1-1}\le i_s\le j_s\le j_{r_1-1}$, for $s>r_1$ such that  $\bos(r_1,s)\in\bs^\circ\sqcup\bs$.
   \end{prop}
   \begin{pf} If  $i_{r_1+1}<i_{r_1-1}$ then since  $\bos$ is connected one of the following holds: $$i_{r_1+1}< i_{r_1-1}<i_{r_1}\le\min\{j_{r_1+1}, j_{r_1-1}\}\ \ {\rm or}
\ \ i_{r_1}<i_{r_1+1}<i_{r_1-1}\le j_{r_1}<\min\{j_{r_1+1}, j_{r_1-1}\}.$$ Since $[i_{r_1+1}, j_{r_1+1}]$ and $[i_{r_1-1}, j_{r_1-1}]$ do not overlap we get $j_{r_1-1}\le j_{r_1+1}$ which contradicts the assumption that $\bos$ is  stable and the proposition is proved for $s=r_1+1$. 
Assume that $\bos(r_1,s)\in \bs^\circ\sqcup\bs$ and that we have proved the result for $s-1$ with $\bos(r_1,s-1)\in \bs^\circ\sqcup\bs$. Then either $i_{s-1}<i_s$ or $j_s<j_{s-1}$ and hence, using that $\bos$ is connected,  one of the following holds: $$i_{r_1-1}\le i_{s-1}<i_s\le j_{s-1}\leq j_{r_1-1}\ \ \ {\rm or}\ \ \ i_{r_1-1}\le i_{s-1}\le j_s<j_{s-1}\leq j_{r_1-1}.$$ Since $[i_{r_1-1}, j_{r_1-1}]$ and $[i_s, j_s]$ do not overlap we get $j_s\le j_{r_1-1}$ in the first case and $i_{r_1-1}\le i_s $ in the second case and the proof of the proposition is complete.
        \end{pf}
       \begin{cor}
            Assume that $\bos$ is connected and stable and suppose that $2\le p\le r-1$ is such that $\bos(p-2,p+1)\notin\bs^\circ\sqcup\bs$.
            Then $i_{p-1}\le i_s<j_s\le j_{p-1}$ for all $s>p$ such that $\bos(p,s)\in\bs^\circ\sqcup\bs$.  In particular if $\bos(0,2)\in\bs^\circ $ then $j_1\ge  j_s$ for all $s\ge 2$  with strict inequality holding if $r_1>2$.
        \end{cor}
        \begin{pf} The first assertion of the corollary follows by working with $\bos(p-2, r)$.  For the second one we note that $j_1\ge j_s$ for all $2\le s\le r_1$ and that $j_{r_1-1}\ge j_s$ if $\bos(r_1, s)\in\bs$. If $r_2$ is the maximal value of $s$ with this property then $\bos(r_2-1, r_2+1)\in\bs^\circ$ and hence by induction we get $j_{r_2}\ge j_\ell$ for all $\ell>r_2$. Iterating it follows that $j_1\ge j_s$ for all $s\ge 2$. If $r_1>2$ then we have $j_1>j_{r_1-1}\ge j_{r_2}\ge j_\ell$ and the proof is complete.
        \end{pf}

\subsection{} \begin{prop}\label{spdj} Suppose that $\bos\in \uS$ is stable and that $\bos(0,r_1)\in\bs^\circ$. For $1\leq \ell<p\leq r_1<s\leq r$ the following pairs of intervals do not overlap:
\begin{gather} \label{ip}([i_p,j_1], [i_\ell, j_{\ell+1}]),\ \ 1\le \ell<p, \ \ \  ([i_p,j_1], [i_s,j_s]),\ \ s>r_1,\\
  \label{is} ([i_s, j_s], [i_\ell, j_{\ell+1}]),\  1\le \ell<\min\{p,r_1-1\}\ \  {\rm and}\ \ s>r_1,\\ \label{is1} ([i_s, j_s], [i_{r_1-1}, j_{r_1}]),\ \ \bos(r_1-1,s)\notin\bs.\end{gather}
\end{prop}
\begin{pf} If $1\le \ell<p$ then $i_p<i_\ell\le j_{\ell+1}<j_1$ showing that the first pair of intervals in \eqref{ip} do not overlap. 
If $p<r_1$ and $s>r_1$  or if $p=r_1$ and $\bos(r_1-1,s)\notin\bs^\circ\sqcup\bs$ then Lemma \ref{tauprops} proves that $[i_p,j_1]$ and $[i_s, j_s]$ do not overlap if $s>r_1$.
If  $p=r_1$ and $\bos(r_1-1,s)\in\bs^\circ\sqcup\bs$ then  Proposition \ref{structure} gives  $i_{r_1}<i_{r_1-1}\le i_s<j_s\le j_{r_1-1}\le j_1$. This completes the proof that the intervals in \eqref{ip} do not overlap.\\\\
The fact that the intervals in \eqref{is} do not overlap  is immediate since $\ell+1<r_1$ and hence $[i_s,j_s]$, $[i_\ell, j_\ell]$ and $[i_{\ell+1}, j_{\ell+1}]$ satisfy the conditions of Lemma \ref{tauprops}.\\\\
Finally if $\bos(r_1-1, s)\notin\bs$ then $[i_s,j_s]$ does not overlap $[i_{r_1-1}, j_{r_1-1}]$ and $[i_{r_1}, j_{r_1}]$ and an application of Lemma \ref{tauprops}
shows that the intervals in \eqref{is1} do not overlap.
 \end{pf}
 
 \subsection{} We record the following for later use.
 \begin{lem}\label{tildesp} Let $\bos\in\uS^{\pr}$ be stable with  $\bos(0,r_1)\in\bs^{\circ}$. For $1 \le  p<r_1$ we have that $\tilde\bos_p=([i_p,j_1])\vee \bos(p,r)\in \uS$ and is connected. Moreover  $\tilde\bos_p(0,2)$ is contained in a prime factor of $\tilde\bos_p$.\end{lem}
  \begin{pf}  If $p=1$ then $\tilde\bos_1 = \bos$ and there is nothing to prove. If $1<p<r_1$ it follows from Corollary  \ref{structure}  that $j_1>j_s$ for all $s\ge 2$. Hence $\tilde\bos_p$ satisfies the first condition in the definition of $\uS$.
  The second condition holds since  $p<r_1$ and so $i_{p+1}<i_p\le j_{p+1}<j_1$ and hence $\tilde\bos_p$ is connected by Lemma \ref{elemalt}(i).
  Finally, \eqref{ip} shows that the third condition is also satisfied. \\\\
Suppose that $\tilde\bos_p(0,2)$ is not contained in a prime factor of $\tilde\bos_p$. Since $\tilde\bos_p$ is connected, an inspection of Definition \ref{pfdef} shows that we must have $ p+1=r_1$ and $i_p=i_{p+2}$ or $j_1=j_{p+2}$. Since $\bos$ is prime the first cannot happen and the second fails  since $j_1>j_{p+2}$ by Corollary \ref{structure}.  
\end{pf}

\subsection{Proof of Proposition \ref{sps}} Recall that $\bos\in\uS^{\pr}$ is stable and that we have to prove that $\bos_p=([i_1,j_2], [i_2,j_3],\cdots ,[i_{p-1}, j_p])\vee \bos(p,r)$ is a stable alternating snake. Since $\bos_1=\bos(1,r)$ the result is immediate from Lemma \ref{elemalt}(i) when $p=1$. From now on we assume that $p\geq 2$ and that $\bos(0,r_1)\in\bs^\circ$. The case $\bos(0,r_1)\in\bs$ follows by working with $\Omega(\bos)$ and using Lemma \ref{elemalt}(ii). \\\\
 To show that $\bos_p\in\uS$ it suffices to prove the following three statements:
\begin{enumerit}
    \item[(i)] For $1\le \ell<p<s$ either $i_\ell\ne i_s$ or $j_{\ell+1}\ne j_s$.
    \item[(ii)] $i_{\ell}<i_{\ell-1}$ and $j_{\ell+1}<j_{\ell}$ 
 if $1< \ell<p$;   $i_{p+1}<i_{p-1}$ and $j_{p+1}<j_{p}$ if $p<r_1$; and $i_{r_1-1}<i_{r_1+1}$ and $j_{r_1}<j_{r_1+1}$ if $p=r_1$.
    \item[(iii)] If $\bos_p(\ell-1, s)\notin\bs^\circ\sqcup\bs$ for some $1\le \ell <p<s $ then $[i_\ell, j_{\ell+1}]$ and $[i_s,j_s]$ do not overlap.\\\\
\end{enumerit}
 It  suffices to prove (i) when $\ell=1$; working with $\bos(\ell-1,r)$ then gives the result for $2\le \ell \le r_1$. Notice that $j_s< j_2$ if $p<s\le r_1$ and hence it suffices to prove (i) when $s>r_1$.
  If $r_1>2$ then  $j_{r_1+1}<j_{r_1-1}\le j_2$ and hence we can consider $s>r_1+1$. Proposition \ref{spdj} (equations \eqref{is} and $\eqref{is1}$) gives that the intervals $[i_{1},j_{2}]$ and $[i_{s-1},j_{s-1}]$ do not overlap. Since $[i_s,j_s]$ and  $[i_{s-1},j_{s-1}]$ do overlap it follows that if $i_1=i_s$ then $j_2\ne j_s$ and vice versa.
 If $r_1=2$ then  $j_2<j_{s}$ for all $s$ with $\bos(1,s)\in\bs$ and hence we may assume that $s$ is such that $\bos(1,s)\notin\bs$.   If in addition $\bos(1,s-1)\notin \bs$ then arguing as in the $r_1>2$ case we have that  $[i_1,j_2]$ and $[i_{s-1},j_{s-1}]$ do not overlap and hence either $i_1\ne i_s$ or $j_2\ne j_s$. Hence it remains to consider the case when $r_1=2$ and $\bos(1,s-1)\in\bs$ with $\bos(1,s)\notin\bs$; we claim that $s=4$. In fact, if $s>4$ with $i_s=i_1$ and $j_s=j_2$ we would have 
 $$i_s=i_1<i_3\leq j_2=j_s<j_3,$$
 contradicting the fact that $[i_s,j_s]$ and $[i_3,j_3]$ do not overlap since $\bos(2,s)\notin\bs^\circ\sqcup\bs$, proving the claim. Finally, noting that $j_4\ne j_2$, since $\bos$ is prime, the proof of (i) is complete.
 \\\\
If $\ell<p$ or if $p<r_1$ then part (ii) holds since $\bos(0,r_1)\in\bs^\circ$ while if $p=r_1$ the assertion holds by Proposition \ref{structure} and \eqref{ijneq}, since $\bos$ is prime and stable.\\\\
It suffices to prove part (iii) when $\ell=1$; working with $\bos(\ell-1,r)$ then gives the result for $2\le \ell \le r_1$.
Note that part (ii) implies that $\bos_p(0, s)\in\bs^\circ\sqcup\bs$
if $s<r_1$.  Hence we may assume that $s\ge r_1$ in which case we have to prove that $[i_1,j_2]$ and $[i_{s+1}, j_{s+1}]$ do not overlap.
 If $r_1>2$ this  follows from Lemma  \ref{tauprops} applied to $[i_1,j_1], [i_2,j_2], [i_{s+1},j_{s+1}]$. If $r_1=2$ then $\bos_p(0,s)=([i_1, j_2],[i_3,j_3],\cdots, [i_{s+1}, j_{s+1}])\notin\bs^\circ\sqcup\bs$ only if $\bos(1,s+1)\notin\bs$ and hence the result again follows from Lemma  \ref{tauprops} applied to $[i_1,j_1], [i_2,j_2], [i_{s+1},j_{s+1}]$.\\\\
 Finally we prove that $\bos_p$ is stable. If $1\le p<r_1-1$ then $$\bos_p(0,r_1-1)=([i_1, j_2], \cdots, [i_{p-1}, j_p], [i_{p+1}, j_{p+1}],\cdots, [i_{r_1-1}, j_{r_1-1}], [i_{r_1},j_{r_1}])$$ and the result follows since $\bos$ is stable. If $p=r_1-1$ we have 
 $$\bos_p(0,r_1)=([i_1, j_2], \cdots, [i_{r_1-2}, j_{r_1-1}], [i_{r_1},j_{r_1}], [i_{r_1+1},j_{r_1+1}]). $$
Since  $[i_{r_1+1}, j_{r_1+1}]$ and $[i_{r_1-2}, j_{r_1-2}]$ do not  overlap and $j_{r_1+1}<j_{r_1-1}<j_{r_1-2}$, by Proposition \ref{structure}, it follows that if   $i_{r_1+1}<i_{r_1-2} $  then we must have that $j_{r_1+1}<i_{r_1-2}$.\\\\
Finally if $p=r_1$ then
$$\bos_p(0,r_1)=([i_1,j_1],\cdots,[i_{r_1-2}, j_{r_1-1}], [i_{r_1-1}, j_{r_1}], [i_{r_1+1}, j_{r_1+1}]).$$ Here we have $j_{r_1+1}<j_{r_1-1}$, by Proposition \ref{structure}, and we  must check that $i_{r_1+1}<i_{r_1-2}$ forces $j_{r_1+1}<i_{r_1-2}$. In addition if  $\bos(r_1-1, r_1+2)\notin\bs$ we have $j_{r_1+2}<j_{r_1}$ and   we must check that $i_{r_1+2}<i_{r_1-1}$ forces $j_{r_1+2}<i_{r_1-1}$. The  assertions follow from the fact that $[i_{r_1-2}, j_{r_1-2}]$ and $[i_{r_1+1}, j_{r_1+1}]$  do not overlap, and  if  $\bos(r_1-1, r_1+2)\notin\bs$  then $[i_{r_1}, j_{r_1}]$ and $[i_{r_1+2}, j_{r_1+2}]$
do not overlap and the following inequalities$$i_{r_1+1}<i_{r_1-2}\le j_{r_1-1}<j_{r_1-2}\ \ {\rm and}\ \ i_{r_1+2}<i_{r_1-1}\le j_{r_1}<j_{r_1-1}. $$ The proof of the proposition is complete.
\\\\

\subsection{}\label{spnotconn} Suppose that  $\bos\in\uS^{\pr}$ and that  $\bos(0,r_1)\in\bs^\circ$. It   is clear that $\bos_p$ is connected if and only if $([i_1,j_2],\cdots, [i_{p-1}, j_p], [i_{p+1}, j_{p+1}])$ is connected. In turn this is equivalent to the assertion that $\bos_p$ is connected if and only if $i_{m-1}\le j_{m+1}$ for all $2\le m\le p$. \\\\
We need another formulation of this equivalence; namely $\bos_p$ is not connected if and only if there exists $2\le m\le p$ such that $j_s<i_{m-1}$ for all $s\ge m+2$ with $\bos(r_1-1,s)\in\bs$. If $\bos_p$ is not connected then there exists  $2\le m\le p$ minimal such that $j_{m+1}\le i_{m-1}$. Since $\bos(0,r_1)\in\bs^\circ$ it follows that $j_s<i_{m-1}$ for all $m+1\le s\le r_1$. Now using $i_{r_1+1}\le j_{r_1}<\min\{i_{m-1},j_{r_1+1}\}<j_{m-1}$ and the fact that $[i_{m-1},j_{m-1}]$ and $[i_{r_1+1}, j_{r_1+1}]$ do not overlap gives $j_{r_1+1}<i_{m-1}$. Repeating with $i_{r_1+2}\le j_{r_1+1}<\min\{j_{r_1+2},i_{m-1}\}<j_{m-1}$ and further iterations gives the result. The  converse direction is immediate.

 \section{An identity in $\cal K_0(\mathscr F_n)$ and Proof of Proposition \ref{altsumid}}\label{id}
Recall our   convention that $$[V(\bomega_{i,j}\bomega)]=0,\ \ {\rm for\ all}\ \  \bomega\in\cal I_n^+,\ \  {\rm if}\ [i,j]\notin\mathbb I_n.$$
\subsection{} Proposition \ref{altsumid} is immediate form the following stronger result.
\begin{prop}\label{idp}
    Suppose that $\bos\in\uS^{\pr}$ is stable and  $1\le p\le r_1$. Then,
    \begin{enumerit}
        \item[(i)] If  $\bos(0,r_1)\in\bs^\circ$ 
     we have 
     \begin{gather*}
[V(\bomega_{i_p,j_1})][V(\bomega_{\bos_p})] = [V(\bomega_{i_p,j_1}\bomega_{\bos_p})] + 
(1-\delta_{r_1,p})[V(\bomega_{i_{p+1}, j_1}\bomega_{\bos_{p+1}})].\end{gather*} 
\item[(ii)]If $\bos(0,r_1)\in \bs$  we have 
\begin{gather*}
[V(\bomega_{i_1,j_p})][V(\bomega_{\bos_p})] =  [V(\bomega_{i_1,j_p}\bomega_{\bos_p})] + 
(1-\delta_{r_1,p})[V(\bomega_{i_{1}, j_{p+1}}\bomega_{\bos_{p+1}})].\end{gather*}\end{enumerit}
\end{prop}
\begin{pf}
  Lemma \ref{Omega} shows that we can deduce part (ii)  from part (i) by applying $\tilde\Omega$ to both sides of the equality. Hence from  now on we shall assume that $\bos(0,r_1)\in\bs^\circ$.
 If $[i_p,j_1]\notin\mathbb I_n$ then we have $j_1-i_{p+1}>j_1-i_p> n+1$ and the proposition is obviously true. So we further assume from now on that $1\le p\le r_1$ is such that $[i_p,j_1]\in\mathbb I_n$.
\\\\
By Proposition \ref{spdj} the following  pairs of intervals 
$([i_p,j_1], [i_\ell, j_{\ell+1}])$ for $ \ell<p$, $([i_p,j_1], [i_s,j_s])$ for $s>r_1$, and $([i_s, j_s], [i_\ell, j_{\ell+1}])$ for $\ell+1<r_1<s$, do not overlap. Hence 
\eqref{fundtpirr}, Proposition \ref{kkop}(ii) and its corollary give
\begin{gather*}\lie d(V(\bomega_{i_p, j_1}), V(\bomega_{i_\ell, j_{\ell+1}}))=0,\ \ \ell<p,\ \ \lie d( V(\bomega_{i_p,j_1}), V(\bomega_{i_s, j_s}))=0,\ \ s>r_1,\\ \lie d(V(\bomega_{i_s,j_s}),V(\bomega_{i_\ell,j_{\ell+1}}))=0,\ \ \ell+1<r_1<s,\\ \lie d(V(\bomega_{i_p,j_1}), V(\bomega_{\bos_p(0,p-1)}))=0= \lie d(V(\bomega_{i_p,j_1}), V(\bomega_{\bos(r_1,r)})).\end{gather*} 
Corollary \ref{kkop} further gives \begin{gather}\notag\lie d(V(\bomega_{i_{r_1},j_1}), V(\bomega_{\bos_{r_1}}))\le \lie d(V(\bomega_{i_{r_1},j_1}), V(\bomega_{\bos_p(0,r_1-1)}))+ \lie d(V(\bomega_{i_{r_1},j_1}), V(\bomega_{\bos(r_1,r)}))=0,\\ \notag\\ 
\label{ineqsp1}\lie d(V(\bomega_{i_p,j_1}), V(\bomega_{\bos_p}))\le \lie d(V(\bomega_{i_p,j_1}), V(\bomega_{\bos(p,r_1)}))\le 1,\ \ \ p<r_1.
 \end{gather}
 The final inequality in \eqref{ineqsp1} follows from Proposition \ref{naoisnake} once we note that $([i_p,j_1])\vee \bos(p,r_1)\in\bs^\circ$. 
  If  $j_1-i_{p+1}>n+1$ then Proposition \ref{naoisnake} gives 
$\lie d(V(\bomega_{i_p, j_1}),  V(\bomega_{\bos(p,r_1)}))=0 $ and Proposition \ref{idp} follows in this case.
\\\\
 To complete the proof we  consider the cases when  $p<r_1$ and $j_1-i_{p+1}\leq n+1$. 
  The inequalities in \eqref{ineqsp1} and Proposition \ref{kkop} show that the module  $V(\bomega_{i_p,j_1})\otimes V(\bomega_{\bos_p})$ has length at most two. We  prove that it has length exactly two by showing that $V(\bomega_{i_{p+1}, j_1}\bomega_{\bos_{p+1}})$ is a Jordan--Holder component. Noticing that $\bomega_{\bos_p}= \bomega_{\bos_p(0,p-1)}\bomega_{\bos(p,r)}$ and using Lemma \ref{gammaroot} we get
\begin{gather*}\bomega_{i_p,j_1}\bomega_{\bos_p}= \bomega_{i_{p+1}, j_1}\bomega_{\bos_{p+1}}\gamma_0,\ \ 
  {\rm where}\ \ \gamma_0=\prod_{i=i_{p+1}}^{i_p-1}\prod_{j=j_{p+1}}^{j_1-1}\balpha_{i,j}.
  \end{gather*}
   Lemma \ref{tildesp} asserts that $\tilde\bos_p=([i_p,j_1])\vee\bos(p,r)\in\uS $  and  that $([i_p,j_1], [i_{p+1}, j_{p+1}])$ is contained in a prime factor of $\tilde\bos_p$. Hence     Proposition \ref{lwtprime} (with $\bos$  replaced by $\tilde\bos_p$ and $p$ replaced with $1$) gives, $$\bomega_{i_p,j_1}\bomega_{\bos(p,r)}\gamma_0^{-1}\in\wt_\ell(V(\bomega_{i_p,j_1})\otimes V(\bomega_{\bos(p,r)}))\setminus \wt_\ell V(\bomega_{i_p,j_1}\bomega_{\bos(p,r)}).$$
   It follows that $$\bomega:=\bomega_{i_p,j_1}\bomega_{\bos_p}\gamma_0^{-1}\in\wt_\ell^+ M,\ \ M=  V(\bomega_{\bos_p(0,p-1)})\otimes V(\bomega_{i_p,j_1})\otimes V(\bomega_{\bos(p,r)}).$$ Suppose that $\bomega$ is an $\ell$--weight in  some Jordan--Holder component of $M$.  
Then
there  exists   $\gamma\in\cal Q_n^+$  such that $\gamma\preccurlyeq\gamma_0$ and $\bomega_{i_p,j_1}\bomega_{\bos_p}\gamma^{-1}\in\wt_\ell^+ M$.
Write $$\bomega_{i_p,j_1}\bomega_{\bos_p}\gamma^{-1}=\bomega_1\bomega_2\bomega_3,\ \ \bomega_1\in\wt_\ell V(\bomega_{\bos_p(0,p-1)}),\ \ \bomega_2\in\wt_\ell V(\bomega_{i_p,j_1}),\ \ \bomega_3\in\wt_\ell V(\bomega_{\bos(p,r)}),$$
 Now writing $\gamma$  in terms of the generators of $\cal Q_n^+$ we see that $\gamma$ cannot involve any element of the form $\balpha_{i_\ell, j_\ell+1}$ since $i_\ell>i_p-1$ if $1\le \ell\le p-1$. Hence  the discussion in Section \ref{gammcomp} gives $\bomega_1=\bomega_{\bos_p(0,p-1)}$.
 Further if we write $\gamma$ in terms of the generators of $\cal I_n^+$ then  Lemma \ref{gammcomp} shows that $\bomega_{i,j}^{-1}$ can occur in it only if $i_{p+1}\le i\le i_p$ and hence  $\bomega_{\bos_p(0,p-1)}\bomega_{i,j}^{-1}\notin\cal I_n^+$. 
It follows that $$\bomega_{i_p,j_1}\bomega_{\bos(p,r)}\gamma^{-1}\in\wt_\ell^+( V(\bomega_{i_p,j_1})\otimes V(\bomega_{\bos(p,r)})).$$ Using \eqref{gammdom} applied to $\tilde\bos_p=([i_p,j_1])\vee\bos(p,r)$ we get  $\gamma=\gamma_0$. Hence  we have shown that \begin{equation}\label{step1} \bomega\in\wt^+_\ell(V(\bomega_{i_p,j_1})\otimes V(\bomega_{\bos_p})).\end{equation} 
 Proposition \ref{gammacond} applied to $\tilde\bos_p$ gives 
$$ \dim (V(\bomega_{i_p,j_1}\bomega_{\bos(p,r)}))_{\bomegas_{i_{p},j_1}\bomegas_{\bos(p,r)}\gamma_0^{-1}}=0$$  and so the preceding arguments also give  \begin{gather*}\bomega\notin\wt_\ell V(\bomega_{\bos_p(0,p-1)})\otimes V(\bomega_{i_p,j_1}\bomega_{\bos(p,r)})),\\ \gamma\prec\gamma_0,\ \ \bomega_{i_p,j_1}\bomega_{\bos_p}\gamma^{-1}\in\cal I_n^+\implies \gamma_0\gamma^{-1}\notin\cal Q_n^+.\end{gather*}Hence the $\widehat\bu_n$--submodule generated by the  weight space corresponding to $\bomega_{i_p,j_1}\bomega_{\bos}\gamma_0^{-1}$ is irreducible and gives the second Jordan Holder component of 
$V(\bomega_{i_p,j_1})\otimes V(\bomega_{\bos(p,r)})$. This completes the proof of the proposition.  \end{pf}

\section{Proof of Proposition \ref{matrixprop}} Throughout this section we assume that $\bos=([i_1,j_1],\cdots, [i_r, j_r])\in\uS$ is  stable. 
\subsection{} We begin with some preliminary comments. For $1\le p\le r$ we can write $A(\bos)$ as a block matrix where the diagonal blocks are $A(\bos(0,p))$ and $A(\bos(p,r))$, i.e., $$A(\bos)=\begin{bmatrix}
    A(\bos(0,p))& B \\
    C & A(\bos(p,r)) 
\end{bmatrix}.$$ This is clear from the definition if $p=1$ and for $p>1$ a straightforward induction on $r$ gives the result. Moreover if $\bos(p-1,p+1)$ is not connected we have $B=0= C$.\\\\
It is convenient to define $r_1$,  $r_2$ and $r_3$ (if they exist)  to be maximal so that 
$$\bos(0,r_1)\in \bs^\circ\sqcup\bs, \ \ \bos(r_1-1,r_2)\in\bs^\circ\sqcup\bs\ \ {\rm and}\ \ \bos(r_2-1,r_3)\in\bs^\circ\sqcup\bs. $$

 \subsection{Proof of Proposition \ref{matrixprop}(i)} We prove the proposition when $\bos(0,2)\in\bs^\circ$. An application of Lemma \ref{elemalt}(ii) gives the result when $\bos(0,2)\in \bs$. Recall that for this proposition we are assuming also that $\bos$ is prime and hence by Remark \ref{sps} the matrix $A(\bos_p)$ is defined. 
 We prove the proposition  by induction on $p$; since $\bos_1=\bos(1,r)$ it is clear that induction begins at $p=1$.\\\\
 For the inductive step suppose first that $\bos_p$ is not connected. By the discussion in Section \ref{spnotconn} there exists $2\le m\le p$ such that $j_s< i_{m-1}$ for all $m+1\le s\le r_2$. This gives $A(\bos)_{\ell,s}=0$ if $ 1\le \ell\le m-1$ and $s\ge m+1$ and hence $A_p(\bos)$ has a block decomposition \begin{equation}\label{minor1} A_p(\bos)=\begin{bmatrix} A_m(\bos(0,m))&0\\
 C& A_{p-m+1}(\bos(m-1, r)).\end{bmatrix}\end{equation}
On the other hand since $\bos_p(0,m-1)=([i_1,j_2],\cdots, [i_{m-1}, j_m])\in\bs^\circ$ is connected by the  minimality of $m$ and $\bos_p(\ell-1, s)$ is not connected for all $2\le \ell\le m$ and $s\ge m$ we have \begin{equation}\label{minorm}A(\bos_p(0,m-1))= A_m(\bos(0,m)),\ \ A(\bos_p)_{\ell, s}=0= A(\bos_p)_{s,\ell}\ \ 1\le \ell\le m-1,\ \ s\ge m\end{equation}
and so,
\begin{equation}\label{minorsp} A(\bos_p)=\begin{bmatrix} A(\bos_p(0,m-1))&0\\ 0&A(\bos(m-1,r)_{p-m+1}).
    \end{bmatrix}\end{equation} The inductive hypothesis gives 
    $\det A_{p-m+1}(\bos(m-1,r))=\det A(\bos(m-1,r))_{p-m+1})$ and hence the inductive step follows  from \eqref{minor1}, the first equality in \eqref{minorm} and \eqref{minorsp}.\\\\
    We prove the inductive step when $\bos_p$ is connected. The  definition of $A(\bos)$ and $A(\bos_p)$ give \begin{gather*}A(\bos_p)_{1,\ell}=A_p(\bos)_{1,\ell}= \begin{cases} [V(\bomega_{i_1,j_{\ell+1}})],& 1\le \ell\le r_2-1,\\
     0,& {\rm otherwise} ,\end{cases}  \\ A(\bos_p)_{\ell,1}=A_p(\bos)_{\ell, 1}=\begin{cases} [V(\bomega_{i_{\ell+1}, j_2})],& 1\le \ell< r_{1+2\delta_{r_1,2}},\\
     0,& {\rm otherwise}.\end{cases}\end{gather*}
     Since $p\ge 2$ we have 
$\bos_p=([i_1,j_2])\vee\bos(1,r)_{p-1}$ and it is easy to check that\begin{gather} \label{sl2}A_p(\bos)_{s,\ell}=A_{p-1}(\bos(1,r))_{s-1,\ell-1},\  \ A(\bos_p)_{s,\ell}= A(\bos(1,r)_{p-1})_{s-1,\ell-1},\ \ s,\ell\ge 2.
     \end{gather}
     Since $\bos_{p-1}(1,r)$ is also connected we get by the preceding arguments  that the first column and row of $A_{p-1}(\bos(1,r))_{s-1,\ell-1}$ and $A(\bos(1,r)_{p-1})_{s-1,\ell-1}$ are equal. Iterating it follows that $A(\bos_p)=A_p(\bos)$ if $\bos$ is connected and the proof of the proposition is complete.

\subsection{}  

\begin{prop}\label{0entries} 
 Assume that $\bos$ is connected and suppose that $2\le m\le r-1$ be such that $\bos(m-2,m)\in\bs^\circ$ and $\bos(m-2, m+1)\notin\bs^\circ$.
  Then, \begin{gather*}
     \det A(\bos)=\begin{cases} \det A(\bos(0,m-1))\det A(\bos(m-1,r)),& i_{m-1}=i_{m+1},\\
     \det A(\bos(0,m))\det A(\bos(m,r)),& j_{m-1}=j_{m+1}.\end{cases}
     \end{gather*}\end{prop} 
\begin{pf} 
Assume that $i_{m-1}=i_{m+1}$
and let $m^\dagger$ be maximal such that $\bos(m-1,m^\dagger)\in\bs$. 
We first show that $m^\dagger>m+1$; otherwise, using Corollary \ref{structure}, we would  have
$i_{m+2}<i_{m-1}=i_{m+1}\le j_{m+2}\leq j_m<j_{m-1}
$ contradicting the fact that $[i_{m-1}, j_{m-1}]$ and $[i_{m+2}, j_{m+2}]$ do not overlap.  
\\\\
Writing
$$A(\bos)=\begin{bmatrix}A(\bos(0,m-1)) &B\\
C&A(\bos(m-1,r))\end{bmatrix},$$ we claim that 
\begin{itemize}
\item the only non--zero entries in $B$ are in the $(m-1)$-th row and the first $(m^\dagger-m+1)$ columns;
\item  the second row of $C$ is zero.
\end{itemize}
The claim is equivalent to
\begin{gather}\label{0i1} A(\bos)_{m-1,\ell}=[V(\bomega_{i_{m-1},j_\ell})],\ \ m\le \ell\le m^\dagger,\ \ A(\bos)_{m-1,\ell}=0,\ \ m^\dagger<\ell,\\ \label{0i2}
A(\bos)_{s,\ell}= 0,\ \ 1\le s<m-1<\ell, \ \ A(\bos)_{m+1,s}=0,\ \ 1\le s\le m-1.\end{gather} 
Assuming the claim we prove the proposition in this case as follows.
Since $\bos(m-1,m^\dagger)\in\bs$, the  definition of $A(\bos(m-1,r))$ gives 
\begin{gather*}A(\bos(m-1,r))_{2,\ell'}= A(\bos)_{m+1,\ell'+m-1}=[V(\bomega_{i_{m+1},j_{\ell'+m-1}})],\ \ 1\le \ell'\le m^\dagger-m+1,\\ A(\bos(m-1,r))_{2,\ell'}=0,\ \ \ell'>m^\dagger-m+1.\end{gather*} Subtracting the $(m+1)$-th row of $A(\bos)$ from the $(m-1)$-th row we get that $A(\bos)$ is row equivalent to $$A=\begin{bmatrix}A(\bos(0,m-1)) &0\\
C&A(\bos(m-1,r)),\end{bmatrix}$$ and so  $\det A(\bos)=\det A=\det A(\bos(0,m-1))\det A(\bos(m-1,r))$. 
\\\\
We prove that \eqref{0i1}-\eqref{0i2} hold by induction on  $r$.  We show that induction begins when $m=2$. Since $\bos(0,2)\in\bs^\circ$ and we have 
 proved that $2^\dagger >3$ it follows that $\bos(1,4)\in\bs$.  The definition of $A(\bos)$ gives $$A(\bos)_{1,\ell}=[V(\bomega_{i_1,j_\ell})]\ \ {\rm iff}\ \   1\le \ell \le 2^\dagger,\ \ A(\bos)_{\ell, 1}=0,\ \ \ell\ge 3,$$ which shows that  \eqref{0i1}-\eqref{0i2} holds.  Assume that we have proved \eqref{0i1}-\eqref{0i2}  for $r-1$, in particular they  hold for $\bos(0,r-1)$. For the inductive step we  can further assume that $m>2$; in particular this means that  the inductive hypothesis applies to $\bos(1,r)$. Since $A(\bos)_{s,\ell}=A(\bos(1,r))_{s-1,\ell-1}$ if $s,\ell\ge 2$,  the inductive hypothesis gives the result in these cases. Hence  we only have to prove that,
\begin{equation}\label{01} A(\bos)_{1,\ell}=0\ \ {\rm if}\ \ \ell\ge m\ \ {\rm and}\  \ A(\bos)_{m+1,1}=0. \end{equation}
 Recall that  $r_1$ and $r_2$  are  maximal so that $\bos(0,r_1)\in\bs^\circ\sqcup\bs$ and $\bos(r_1-1, r_2)\in\bs^\circ\sqcup\bs$; clearly $m\ge r_1$. 
 If $\bos(0,r_1)\in\bs$ then $m>r_1\ge 2$ since $\bos(m-2,m)\in\bs^\circ$. It follows that $m\ge r_2$ and so the the equalities  in \eqref{01}  hold by the definition of $A(\bos)$. 
 \\\\
 If $\bos(0,r_1)\in\bs^\circ$ then 
  the second equality in \eqref{01} holds by definition. The first also holds if $m>r_2$. Since $\bos(r_1-1, r_2)\in\bs$ we cannot have $m=r_2$.
  Therefore $m=r_1>2$ and, by Corollary \ref{structure}, we have $i_{r_1-1}=i_{r_1+1}< \min\{i_1,j_{r_1+1}\}<j_1$. 
  Since $[i_1,j_1]$ and $[i_{r_1+1},j_{r_1+1}]$ do not overlap we have $j_{r_1}<j_{r_1+1}<j_1$. Assuming that we have proved that $j_{s-1}<i_1$ for $ r_1<s\le r_2$ we use $i_{s}\le j_{s-1}<\min\{j_s, i_1\}<j_1$ and the fact that $[i_1,j_1]$ and $[i_s, j_s]$ do not overlap to conclude that $j_s<i_1$, for all $r_1\leq s\leq r_2$, which proves the first equality in \eqref{01}.\\\\
 Suppose that $j_{m-1}=j_{m+1}$; this time we write
 $$A(\bos)=\begin{bmatrix}A(\bos(0,m)) &B\\
C&A(\bos(m,r))\end{bmatrix}.$$ Let $m^\bullet<m$ be minimal such that $\bos(m^\bullet-1,m)\in\bs^\circ$. 
We claim that 
\begin{itemize}
\item the only possible non--zero entries in $B$ are in the  first column and the last $m-m^\bullet+1$ rows,
\item the $(m-1)$-th column of $A(\bos)$ is zero unless $m^\bullet\le s\le m$. 
\end{itemize}
The claim is equivalent to 
\begin{gather}\label{0j1}
A(\bos)_{s,m+1}=[V(\bomega_{i_{s},j_{m+1}})],\ \ m^\bullet \le s\le m,\ \ 
A(\bos)_{s,m+1}= 0,\ \ s<m^\bullet,\\  \label{0j2}
A(\bos)_{s,\ell}=0, \ \ s<m+1<\ell,\ \ \   A(\bos)_{s,m-1} = 0, \ {\rm for}\  s<m^\bullet \ {\rm or}\  s>m. 
\end{gather}

Assuming the claim the proof of the proposition is then completed as before by subtracting the $(m-1)$-th column from the $(m+1)$-th column of $A(\bos)$ which makes $A(\bos)$ column equivalent to $$A=\begin{bmatrix}A(\bos(0,m)) &0\\
C&A(\bos(m,r))\end{bmatrix}.$$

We prove that \eqref{0j1}--\eqref{0j2} hold by induction on $r$. First, note that $\bos(m, m+2)\in\bs^\circ$ since otherwise we would have
$$i_{m-1}<i_{m+1}<i_{m+2}\le j_{m+1}=j_{m-1}<j_{m+2},$$ 
where the fist inequality follows from Corollary \ref{structure}; this contradicts the fact that $[i_{m-1}, j_{m-1}]$ and $[i_{m+2}, j_{m+2}]$ do not overlap;  We show that induction begins when $m=2$ in which case the first identity in \eqref{0j1} holds by the definition of $A(\bos)$ and the second one  is vacuously true. For the first identity in \eqref{0j2} we have to show that $A(\bos)_{1,\ell}=0=A(\bos)_{2,\ell}$ if $\ell>3$, which is immediate from the fact that $\bos(m,m+2)\in \bs^\circ$ and from the definition of $A(\bos)$ since $\bos(1,3)\in\bs$. For the second one we have to prove that $A(\bos)_{s,1}=0$, $s>m$, which is again immediate from the definition.\\

Assume we have proved \eqref{0j1}--\eqref{0j2} for $r-1$. For the inductive step we can further assume that $m>2$; in particular the inductive step applies to $\bos(1,r)$ and, similarly as in the previous case, we are left to show that 
\begin{gather}
    A(\bos)_{1,m+1}=0=A(\bos)_{1,m-1},\ \ 1<m^\bullet, \ \ \ A(\bos)_{1,\ell} = 0, \ \ m+1<\ell.
\end{gather}
But these are immediate from the definition of $A(\bos)$ using the fact that  $\bos(m-1,m+2)\notin\bs\sqcup\bs^\circ$, since $\bos(m,m+2)\in \bs^\circ$.
\end{pf}

\subsection{Proof of Proposition \ref{matrixprop}(ii)} By Lemma \ref{elemalt}(ii) we can  write $A(\bos)$ as
\begin{equation*}A(\bos)=\begin{bmatrix}
    A(\bos(0, p))& B_p(\bos)\\
    C_p(\bos)& A(\bos( p,r))
\end{bmatrix}\end{equation*}
If  $\bos(p-1, p+1)$ is not connected then  the  definition of $A(\bos)$  gives $B_p(\bos)=C_p(\bos)=0$  and the proposition is clear.\\\\ 
Therefore we can assume that $\bos$ is connected and that $\bos(p-1,p+1)$ is  contained in a prime factor of $\bos$. 
We prove the result when $\bos(p-2,p)\in\bs^\circ$;  the result in the other case follows by working with $A(\Omega(\bos))$ and using equations \eqref{Aomega1} and \eqref{Aomega2} of Lemma \ref{elemalt}. \\\\ Definition \ref{pfdef} and the fact that $\bos$ is stable now give  that one of the following hold: there exists $2\le m\le r-1$ such that \begin{gather*}p=m-1,\ \ \ \ \bos(m-2, m+1)\notin\bs^\circ\sqcup\bs,\ \
i_{m-1}=i_{m+1}\\ 
p=m,\ \ \bos(m-2,m+1)\notin\bs^\circ\sqcup\bs,\ \ j_{m-1}=j_{m+1},
\end{gather*}
Hence  Proposition \ref{0entries} gives
\begin{gather*}
     \det A(\bos)=\begin{cases} \det A(\bos(0,m-1))\det A(\bos(m-1,r)),& i_{m-1}=i_{m+1},\\
     \det A(\bos(0,m))\det A(\bos(m,r)),& j_{m-1}=j_{m+1},
     \end{cases}
 \end{gather*}
 and the proof of the proposition is complete.

\end{document}